\documentclass[12pt,thmsa]{article}
\usepackage{amsmath}
\usepackage{graphicx}
\usepackage{amsfonts}
\usepackage{amssymb}
\usepackage{mathrsfs}

\newtheorem{proposition}{Proposition}
\newtheorem{lemma}{Lemma}
\newtheorem{theorem}{Theorem}
\newtheorem{corollary}{Corollary}
\newtheorem{assumption}{Assumption}
\newtheorem{remark}{Remark}
\begin{document}

\begin{center}
{\Large \textbf{Variable selection in multivariate regression model for spatially dependent data}}

\bigskip

Jean Roland  EBENDE PENDA, St\'ephane BOUKA 

\bigskip

and  Guy Martial  NKIET 

\bigskip

\textsuperscript{}URMI, Universit\'{e} des Sciences et Techniques de Masuku,  Franceville, Gabon.

\bigskip

E-mail adresses : ebendependa@gmail.com,   stephane.bouka@univ-masuku.org,  guymartial.nkiet@univ-masuku.org.

\bigskip
\end{center}

\noindent\textbf{Abstract.}  This paper deals with variable selection in multivariate linear regression model when the data are observations on a spatial domain being a grid of sites in $\mathbb{Z}^d$ with $d\geqslant 2$. We use a criterion that allows to characterize the subset of relevant variables as depending on two parameters, and we propose estimators for these parameters based on spatially dependent observations. We prove the consistency,  under specified assumptions,   of the method thus proposed.  A simulation study made in order to  assess the finite-sample behaviour of the proposed method with comparison to existing ones is presented.
\bigskip



\section{Introduction}
\noindent Variable selection is known as one of the most relevant problems  in statistical modeling using linear regression model. It  has been widely studied in the literature for many years, the earlier works on it  going back  to the  seventies (e.g., \cite{mallows1973}, \cite{schwarz1978}). Traditional   methods for dealing with it are  based on hypothesis testing in forward selection, backward elimination and stepwise selection (see \cite{draper98}), but also   on criteria to be optimized  among subsets of variables  such as Akaike's information criterion (AIC), Schwarz's Bayesian information criterion (BIC) or Mallows $C_p$  (e.g., \cite{fujikoshi1997}). However, more recent methods rely on penalized likelihood or least squares which led to powerful approaches  that are now very popular such as the least absolute shrinkage and selection operator (LASSO) introduced in \cite{tibshirani1996} and   the method of \cite{fan2001} based on penalization of the likelihood by   a symmetric and nonconcave penalty function. Other recent approaches can be found in \cite{an2013}, which proposed  a method of adaptive sparse canonical correlation analysis (ASCCA) obtained by re-casting the multivariate regression problem as a classical canonical correlation analysis  problem for which a  least squares type formulation is  constructed, so that an   adaptive LASSO type penalty together with a  BIC-type selection criterion are applied directly, and in \cite{mbina2023} which introduced a criterion allowing to describe the set of relevant variables as dependent on two parameters which are then estimated, so reducing the variable selection problem to that of estimating this set.  All these  methods have been proposed in order to achieve variable selection in linear regression when the available data are independent  observations. For the case of spatial data, the literature is not very developed despite the interest of this case in some applications such as environmental, ecological and geographic studies where  it is often of importance to identify factors that are associated with a given  index.  However, there exist some works that deal  with variable selection in linear regression model for spatial data. More specifically, dealing with geospatial data,  \cite{hoeting2006} applied AIC   whereas \cite{liao2019} developed a spatial Mallows criterion and \cite{huang2010}  proposed the use of  LASSO to simultaneously select variables, choose
neighborhoods, and estimate parameters. Furthermore, \cite{chu2011} developped a penalized maximum likelihood estimation  and \cite{wang2009} studied penalized least squares procedures, both enabling    simultaneous variable selection and parameter estimation. In this paper, we extend the approach of \cite{mbina2023} to the case of  spatially dependent observations. More specifically, we use a characterization of the subset of relevant variables  by means of a criterion introduced in \cite{mbina2023} and propose an estimation of this subset from a sample observed on a spatial domain being a grid of sites in $\mathbb{Z}^d$ with $d\geqslant 2$, so achieving variable selection. The rest of the paper is organized as follows. In Section \ref{modeli1}, the multivariate linear regression model with spatial observations that is used is specified, the aforementioned criterion is introduced as well as of the resulting characterization of the subset of relevant variables as depending on two parameters: a permutation and a dimension. In Section \ref{laprocedure}, an estimator of the previous criterion based on spatial data is introduced, so leading to estimation of the aforementioned parameters, what achieves the proposed variable selection procedure. Section \ref{res} is devoted to consistency results for our proposal under assumptions    that are  first specified. A simulation study made in order to  assess the finite-sample behaviour of the proposed method with comparison to existing ones is presented in Section \ref{sim}. All the proofs are postponed  to Section \ref{preuves}.

\section{Model and characterization of relevant set}\label{modeli1}

\noindent We  consider the spatial multivariate regression model    defined by:
\begin{equation}\label{model}
	Y _{\mathbf{i}} = BX_{\mathbf{i}} +  \varepsilon_{\mathbf{i}}, \hspace{0.4cm}  \mathbf{i} \in D \subseteq \mathbb{Z}^{d}, \hspace{0.4cm} d\geq2,
\end{equation}
where  $Y _{\mathbf{i}}$ and  $X_{\mathbf{i}}$  are real random vectors valued respectively into  $\mathbb{R}^{q}$ and $\mathbb{R}^{p}$ with $p\geqslant 2$ and $q \geqslant 2)$,  $B$ is a $q\times p$ matrix of  real coefficients,   and $ \varepsilon_{\mathbf{i}}$ is  a centered random vector valued   into $\mathbb{R}^{q}$ and which is  independent of the  $X_{\mathbf{i}}$s. Writing these random vectors and matrix as:	
\[
X_{\mathbf{i}}=\begin{pmatrix} X_{\mathbf{i}}^{(1)} \\\vdots
	\\ X_{\mathbf{i}}^{(p)} \end{pmatrix},\,\,
Y _{\mathbf{i}}=\begin{pmatrix} Y_{\mathbf{i}}^{(1)} \\\vdots
	\\ Y_{\mathbf{i}}^{(q)} \end{pmatrix},\,\,
\varepsilon_{\mathbf{i}}=\begin{pmatrix} \varepsilon_{\mathbf{i}}^{(1)} \\\vdots\\ \varepsilon_{\mathbf{i}}^{(q)} \end{pmatrix},\ \textrm{ and }\  
B=\begin{pmatrix}
	b_{11} & b_{12} &\cdots & b_{1p} \\
	b_{21} & b_{22} &\cdots & b_{2p} \\
	\vdots & \vdots & \ddots & \vdots \\
	b_{q1} & b_{q2} & \cdots & b_{qp}
\end{pmatrix},
\]
it is easily seen that Model   \eqref{model}  is equivalent to 
\begin{equation}\label{model2}
Y _{\mathbf{i}}=\sum_{j=1}^pX_{\mathbf{i}}^{(j)} \textrm{\textbf{b}}_{\bullet j}+\varepsilon_{\mathbf{i}}
\end{equation}
where
\[
\textrm{\textbf{b}}_{\bullet j} = \left(
 \begin{array}{c}
  b_{1j}\\
  b_{2j}\\
  \vdots \\
  b_{qj}
 \end{array}
\right)\in\mathbb{R}^q.
\]
We are interested with the variable selection problem, that is identifying the $X_{\mathbf{i}}^{(j)}$s that   are not relevant in model \eqref{model2} among $X_{\mathbf{i}}^{(1)},\cdots,X_{\mathbf{i}}^{(p)}$. An  usual approach consists in considering that these variables correspond to   vectors  $\textrm{\textbf{b}}_{\bullet j}$ equal to  $0$. So, putting  $I = \left\{1, \cdots , p\right\}$,  and denoting by $\Vert\cdot\Vert_{\mathbb{R}^k}$  the usual Euclidean norm of $\mathbb{R}^k$, we consider the subset  $I_{0} = \left\{j \in I\,/\, \|\textrm{\textbf{b}}_{\bullet j}\|_{\mathbb{R}^q}= 0  \right\}$ which is assumed to be non-empty, and
we tackle the variable selection problem in Model  \eqref{model}  as a problem of  estimating the subset   $I_1=I-I_0$. Such an approach was considered in \cite{mbina2023}, but on the basis of an i.i.d. sample of the random vectors involved in the regression model. Here, we consider instead spatially dependent data furnished by observations on a spatial domain which is the grid  $\mathcal{I}_{\mathbf{n}}=\{1,\cdots,n\}^d$ in $\mathbb{Z}^d$, with $n\in\mathbb{N}^\ast$. We therefore seek to estimate $I_1$  on the basis of the sample  $\{(X_{\mathbf{i}},Y_{\mathbf{i}})\}_{\mathbf{i}\in\mathcal{I}_{\mathbf{n}}}$. For doing that, we will use a characterization of $I_1$ thanks to a criterion introduced in \cite{mbina2023} and defined as follows. Assuming that each   $(X_{\mathbf{i}},Y_{\mathbf{i}})$ has  the same distribution than a pair $(X,Y)$ of random vectors satisfying  $\mathbb{E}\left(\Vert X\Vert_{\mathbb{R}^p}^2\right)<+\infty$,  $\mathbb{E}\left(\Vert Y\Vert_{\mathbb{R}^q}^2\right)<+\infty$ and, without loss of generality, $\mathbb{E}\left(X\right)=0$,  $\mathbb{E}\left(Y\right)=0$, we consider the covariance and cross-covariance operators given by
\[
V_1=\mathbb{E}\left(X\otimes X\right)\,\,\,\textrm{ and }\,\,\,V_{12}=\mathbb{E}\left(Y\otimes X\right),
\]
where $\otimes$ denotes the usual tensor product such that $u\otimes v$ is the linear map defined by $(u\otimes v)(h)=\langle u,h\rangle v$, and $T^\ast$ denotes the adjoint operator of $T$, and we assume that  $V_1$ is an invertible operator. For any subset $K:=\{k_1,\cdots,k_r\}$ of $I$ (with $r\leqslant p$), we consider the canonical projector $A_K$ from $\mathbb{R}^p$ to $\mathbb{R}^r$ defined as :
\[
A_K\,:\,x=\begin{pmatrix} x_1 \\\vdots
	\\ x_p \end{pmatrix}\in\mathbb{R}^p\,\,
\longmapsto  \begin{pmatrix} x_{k_1} \\\vdots
	\\ x_{k_r} \end{pmatrix}\in\mathbb{R}^r,
\]
and we put
\[
\Pi_{K} = A_{K}^{*}(A_{K}V_{1}A_{K}^{*})^{-1}A_{K},
\]
where $T^{-1}$ denotes the inverse of the operator $T$. Then we define the aforementioned criterion as 
\begin{equation}\label{crit}
\xi_{K}=\left\|V_{12} - V_{1}\Pi_{K}V_{12}\right\|_{\mathcal{H}},
\end{equation}
where $\Vert\cdot\Vert_{\mathcal{H}}$ denotes the norm of operators given by $\Vert T\Vert_{\mathcal{H}}=\sqrt{\textrm{tr}\left(TT^\ast\right)}$. This criterion just is a distance between the matrix $B$ of coefficients  of  \eqref{model} and the matrix $\widetilde{B}_K$  of coefficients    corresponding  to the case where only the variables  $X_{\mathbf{i}}^{(j)}$ such that $j\in K$ are considered in \eqref{model2} (see \cite{mbina2023}). It therefore measures the loss occasionned by selecting the variables whose indices belong to $K$. Thanks to this criterion   $I_1$ can be expressed  as a function of two parameters, a permutation and a dimension, which then have to be estimated in order  to obtain an estimator of  $I_1$. Indeed, it is shown in \cite{mbina2023} that $\xi_K=0$ if and only if $I_1\subset K$. Consequently, for any $i\in I$, putting $K_i=I-\{i\}$ we have $i\in I_1$ if and only if  $\xi_{K_i}>0$. Hence,  sorting  the  $\xi_{K_i}$s in decreasing order, we can identify the elements of $I_1$ as the indices $j$ such that  $\xi_{K_j}>0$. More specifically,  let $\tau$ be the permutation of $I$, that is a one-to-one map from $I$ to itself, satisfying $
\xi_{K_{\tau(1)}}\geqslant\xi_{K_{\tau(2)}}\geqslant\cdots\geqslant\xi_{K_{\tau(p)}}
$
with  $
\tau(i)<\tau(j)$   if $\xi_{K_{\tau(i)}}=\xi_{K_{\tau(j)}}$ and  $i<j$. Since $I_0$ is not empty, there exists an integer  $s$ in $\{1,\cdots,p-1\}$, called dimension, such that 
\[
\xi_{K_{\tau(1)}}\geqslant\xi_{K_{\tau(2)}}\geqslant\cdots\geqslant\xi_{K_{\tau(s)}}>0=\xi_{K_{\tau(s+1)}}=\cdots\xi_{K_{\tau(p)}},
\]
and, therefore,
\begin{equation}\label{i1}
I_1=\{\tau(i)\,/\,1\leqslant i\leqslant s\}.
\end{equation}
Equation \eqref{i1} furnishes an expression of $I_1$ as a function of the permutation $\tau$ and the dimension $s$, and shows that a plug-in  estimate of this set can be obtained from estimators of these two parameters.

 \section{The variable selection procedure}\label{laprocedure}
 
\noindent Our variable selection procedure consists in estimating $\tau$ and $s$ from the sample $\{(X_{\mathbf{i}},Y_{\mathbf{i}})\}_{\mathbf{i}\in\mathcal{I}_{\mathbf{n}}}$ observed on the spatial domain $\mathcal{I}_{\mathbf{n}}$. For doing that, we first need to estimate the criterion given in \eqref{crit}. In this section, estimation of this criterion is adressed, and we then define our proposal for estimating $I_1$ via estimation of $\tau$ and $s$. 
\subsection{Estimation of the criterion}
\noindent We first have to estimate the covariance and cross-covariance operators. Considering the empirical means
$$
\overline{X}^{(\mathbf{n})}=\dfrac{1}{n^{d}}\displaystyle{\sum_{\mathbf{i} \in \mathcal{I}_{\mathbf{n}}}^{} X_{\mathbf{i}}}\,\,\textrm{ and }\,\,
\overline{Y }^{(\mathbf{n})}=\dfrac{1}{n^{d}}\displaystyle{\sum_{\mathbf{i} \in \mathcal{I}_{\mathbf{n}}}^{} Y _{\mathbf{i}}},
$$ 
we estimate $V_1$ and $V_{12}$  by the random operators given, respectively, by
$$\widehat{V}_{1}^{(\mathbf{n})}=\dfrac{1}{n^{d}}\displaystyle{\sum_{\mathbf{i} \in \mathcal{I}_{\mathbf{n}}}^{} (X_{\mathbf{i}}-\overline{X}^{(\mathbf{n})} )\otimes}(X_{\mathbf{i}}-\overline{X}^{(\mathbf{n})})$$  and $$\widehat{V}_{12}^{(\mathbf{n})}=\dfrac{1}{n^{d}}\displaystyle{\sum_{\mathbf{i} \in \mathcal{I}_{\mathbf{n}}}^{} (Y _{\mathbf{i}}-\overline{Y }^{(\mathbf{n})})\otimes (X_{\mathbf{i}}-\overline{X}^{(\mathbf{n})} )}.
$$ 
We will se later that  $\widehat{V}_{1}^{(\mathbf{n})}$ converges almost surely to $V_1$ as $n\rightarrow +\infty$. So, for $n$ large enough, it is also an invertible operator and it is, therefore, possible to put for any subset $K$ of $I$:
$$
\widehat{\Pi}_{K}^{(\mathbf{n})} = A_{K}^{*}(A_{K}\widehat{V}_{1}^{(\mathbf{n})}A_{K}^{*})^{-1}A_{K}.
$$
Then, a plug-in estimator of $\xi_K$ is given by:
$$\widehat{\xi}_{K}^{\mathbf{(n)}} = \left\Vert \widehat{V}_{12}^{\mathbf{(n)}} - \widehat{V}_{1}^{\mathbf{(n)}}\widehat{\Pi}_{K}^{\mathbf{(n)}}\widehat{V}_{12}^{\mathbf{(n)}}\right\Vert_{\mathcal{H}}.$$
\subsection{Estimation of the permutation}
\noindent It would be natural to estimate $\tau$ by sorting the $\widehat{\xi}_{K_i}^{\mathbf{(n)}}$s in decreasing order as it was done for defining $\tau$ from the $\xi_{K_i}$s , but  such an approach does not  guarantee the consistency of the resulting estimator  because  of possible ties in the values of the   $\widehat{\xi}_{K_i}^{(\mathbf{n})}$s. That  is why we rather use estimates of the $\xi_{K_i}$s obtained from appropriate penalizations of the $\widehat{\xi}_{K_i}^{(\mathbf{n})}$s which allow to avoid ex-aequos.   More specifically, we consider
\begin{equation}\label{pen1}
\widehat{\phi}_{i}^{\left(  \mathbf{n}\right)  }=\widehat{\xi}_{K_{i}}^{\left(
\mathbf{n}\right)  }+\frac{f\left(  i\right)}{n^{d\gamma}}, 
\end{equation}
 where
$0<\gamma<1/2  $ and $f$ is  a strictly decreasing function from
$I$ to $\mathbb{R}_{+}$, and we estimate $\tau$ by sorting  the $\widehat{\phi}_{i}^{\left(  \mathbf{n}\right)  }$s in decreasing order. The resulting estimator is the permutation $\widehat{\tau}_\mathbf{n}$ of $I$ satisfying: 
\[
\widehat{\phi}_{\widehat{\tau}_\mathbf{n}(1)}^{\left(  \mathbf{n}\right)  }>\widehat{\phi}_{\widehat{\tau}_\mathbf{n}(2)}^{\left(  \mathbf{n}\right)  }>\cdots>\widehat{\phi}_{\widehat{\tau}_\mathbf{n}(p)}^{\left(  \mathbf{n}\right)  }.
\]
\subsection{Estimation of the dimension}

\noindent For any $i\in I$,  we consider  the subset  $J_i=\{\tau(\ell)\,/\,1\leqslant\ell\leqslant i\}$ of $I$. Then,  we see from \eqref{i1} that $J_i\subset I_1$ if  $i\geqslant s$. Hence   $\xi_{ J_i}=0$ if $i\geqslant s$, and  $\xi_{ J_i}>0$ if $i< s$. Consequently, the dimension $s$  equals  smallest integer $i\in I$ for which $\xi_{ J_i}$ has its minimum value. So, for estimating $s$ we could minimize   $\widehat{\xi}_{\widehat{J}_i}^{(\mathbf{n})}$, where $\widehat{J}_i$ is an estimate of $J_i$, but for the same reason as above, we will use instead a penalized estimator. More precisely, putting     $\widehat{J}_{i}^{ (\mathbf{n} )}=\{\widehat{\tau}_\mathbf{n}(\ell)\,/\,1\leqslant\ell\leqslant i\}$, we consider 
\begin{equation}\label{pen2}
\widehat{\psi}_{i}^{\left( \mathbf{n}\right)  }=\widehat{\xi}_{\widehat{J}_{i}^{ (\mathbf{n} )}}^{(\mathbf{n})}+\frac{g\left(  \widehat{\tau}_\mathbf{n}(i) \right)}{n^{d\beta}},
\end{equation}
where  $0<\beta<1/2  $ and $g$ is  a strictly increasing function  from $I$ to \ $\mathbb{R}_{+}$.
Then, we take as estimator of $s$ the statistic:
\begin{equation*}
\widehat{s}_\mathbf{n}=\arg\min_{i\in I}\left(  \widehat{\psi}_{i}^{\left( \mathbf{n}\right)  }\right).
\end{equation*}

\noindent Finally, having estimated the parameters $\tau$ and $s$, we propose to select relevant variables in Model \eqref{model} by estimating $I_1$ by the set
\[
\widehat{I}_{1}^{\left(  \mathbf{n}\right)  }=\{\widehat{\tau}_\mathbf{n}(i)\,/\,1\leqslant i\leqslant \widehat{s}_\mathbf{n}\}.
\]
\section{Consistency results}\label{res}

\noindent In this section, we first introduce the assumptions needed to obtain the main results of the paper and then we state
convergence theorems that give consistency for our variable selection method.
\subsection{Assumptions}\label{as}

\bigskip

\begin{assumption}\label{ass0}
The process $\{\varepsilon_\mathbf{i}\}_{\mathbf{i}\in\mathbb{Z}^d}$ is centered, independent of  $\{X_\mathbf{i}\}_{\mathbf{i}\in\mathbb{Z}^d}$ and satisfies  $\sup_{\mathbf{i}\in\mathbb{Z}^d}\mathbb{E}\left(\Vert \varepsilon_\mathbf{i}\Vert_{\mathbb{R}^q}^2\right)<+\infty$.  
\end{assumption}

\begin{assumption}\label{ass1}
The process $\{Z_\mathbf{i}\}_{\mathbf{i}\in\mathbb{Z}^d}$, where  $Z_\mathbf{i}=(X_\mathbf{i},Y_\mathbf{i})$, is strictly stationnary, satisfies   $\mathbb{E}\Big(\Vert Z_\mathbf{i}\Vert_{\mathbb{R}^{p+q}}^8\Big)<+\infty$  and  is  strongly $\alpha$-mixing: there exists a decreasing  function $\varphi$ satisfying  $\lim_{t\rightarrow +\infty}\varphi(t)= 0$, such that  
\begin{eqnarray*}
	\alpha(\mathscr{B}(\textbf{S}),\mathscr{B}(\textbf{S}^\prime))&=&\sup\{|\mathrm{P(A\cap B)}-\mathrm{P}(A)\mathrm{P}(B)|, A\in \mathscr{B}(\textbf{S}), B\in \mathscr{B}(\textbf{S}^\prime)\}\\
	&\leqslant &\varphi(\delta(\textbf{S},\textbf{S}^\prime)),
\end{eqnarray*}
where $\textbf{S}$ and $\textbf{S}^\prime$ are finite subsets of $\mathbb{Z}^d$, $\mathscr{B}(\textbf{S})=\sigma(Z_\mathbf{i},\,\mathbf{i}\in\textbf{S})$,  $\mathscr{B}(\textbf{S}^\prime)=\sigma(Z_\mathbf{i},\,\mathbf{i}\in\textbf{S}^\prime)$, and $\delta$ is the distance between finite  subsets of $\mathbb{Z}^d$ defined by $\delta(E,F)=\min\{\Vert \mathbf{i}-\mathbf{j}\Vert_2,\,\mathbf{i}\in E,\,\mathbf{j}\in F\}$ with $\Vert \mathbf{i}-\mathbf{j}\Vert_2=\left(\sum_{k=1}^d(i_k-j_k)^2\right)^{1/2}$.
\end{assumption}

\bigskip

\noindent We say that a strictly stationnary process $\{U_\mathbf{i}\}_{\mathbf{i}\in\mathbb{Z}^d}$ valued into an Euclidean space $E$ with dimension $m$ and inner product $\langle\cdot,\cdot\rangle_E$  is isotropic if,  letting $\{u_k\}_{1\leqslant k\leqslant m}$ be an orthonormal basis of $E$ such that $u_k$ is an eigenvector of the covariance operator 
\[
V_U=\mathbb{E}\Big(\left(U_\mathbf{i} -\mathbb{E}(U_\mathbf{i})\right)\otimes\left(U_\mathbf{i}-\mathbb{E}(U_\mathbf{i})\right)\Big)
\]
associated to  the $k$-th largest eigenvalue $\lambda_k=Var\left(\left\langle U_\mathbf{i} -\mathbb{E}(U_\mathbf{i}),u_k\right\rangle_E\right)$ , we have
\begin{equation*}\label{isotropic}
	\mathbb{E}\Big[\big<U_{\mathbf{i}}-\mathbb{E}(U_\mathbf{i}), u_k \big>_{E }\big<U_{\mathbf{j}}-\mathbb{E}(U_\mathbf{j}), u_{\ell} \big>_{E }\Big] = \rho_{k\ell}\lambda_{k}^{1/2} \lambda_{\ell}^{1/2}\Psi_{k,\ell}(\lVert\mathbf{i}-\mathbf{j} \lVert_2),
\end{equation*}
where the functions  $\Psi_{k,\ell} \colon\mathbb{R}_+\to \mathbb{R}_+$, called spatial correlation functions,  satisfies $\Psi_{k,k} (0)=1$ and  $\sum\limits_{t=1}^{+\infty}t^{d-1}\Psi_{k,\ell}(t)<+\infty$, and $\rho_{k\ell}$ is  a colocated correlation coefficient equal to $1$ if $k=\ell$.

\bigskip

\begin{assumption}\label{ass2}
The process $\{X_\mathbf{i}\}_{\mathbf{i}\in\mathbb{Z}^d}$ is isotropic.
\end{assumption}

\begin{assumption}\label{ass3}
For any $L>0$:
\begin{itemize}
\item[$(i)$]The processes  $\{Z_\mathbf{i}^{(L)}\}_{\mathbf{i}\in\mathbb{Z}^d}$ and  $\{Z_\mathbf{i}^{(L\ast)}\}_{\mathbf{i}\in\mathbb{Z}^d}$, where $Z_\mathbf{i}^{(L)} = Z_\mathbf{i}\textrm{\large \textbf{1}}_{\{ \lVert  Z_\mathbf{i}\lVert_{\mathbb{R}^{p+q}} \leqslant L \}} $ and    	 $Z_\mathbf{i}^{(L\ast)}= Z_\mathbf{i}\textrm{\large \textbf{1}}_{\{ \lVert  Z_\mathbf{i}\lVert_{\mathbb{R}^{p+q}} > L \}} $, are isotropic, the first one having spatial correlation functions $\Upsilon_{k\ell}^{(L)}$ and  colocated correlation coefficients $\upsilon_{k\ell}^{(L)}$  satisfying 
\begin{equation}\label{cvcor}
\lim_{n\rightarrow +\infty}\dfrac{1}{n^{d}} \sum\limits_{\stackrel{\mathbf{i,j} \in \mathcal{I}_{\mathbf{n}}}{\mathbf{i}\neq \mathbf{j}}}\Upsilon_{k\ell}^{(L)}(\lVert\mathbf{i}-\mathbf{j} \lVert_{2})=\zeta_{k\ell}^{(L)},\,\,\,\lim_{L\rightarrow +\infty}\left(\zeta_{k\ell}^{(L)}\right)=\zeta_{k\ell}     \,\,\,\textrm{ and }\lim_{L\rightarrow +\infty}\left(\upsilon_{k\ell}^{(L)}\right)=\upsilon_{k\ell},
\end{equation}
where $\zeta_{k\ell}^{(L)}$, $\zeta_{k\ell}$  and $\upsilon_{k\ell}$ are real numbers.
\item[$(ii)$]The processes  $\{Z_\mathbf{i}^{(L)}\otimes Z_\mathbf{i}^{(L)} \}_{\mathbf{i}\in\mathbb{Z}^d}$ and  $\{Z_\mathbf{i}^{(L\ast)}\otimes Z_\mathbf{i}^{(L\ast)}\}_{\mathbf{i}\in\mathbb{Z}^d}$ are isotropic, the first one having spatial correlation functions $\Theta_{k\ell}^{(L)}$ and  colocated correlation coefficients $\eta_{k\ell}^{(L)}$  satisfying 
\begin{equation}\label{cvcor2}
\lim_{\mathbf{n}\rightarrow +\infty}\dfrac{1}{n^{d}} \sum\limits_{\stackrel{\mathbf{i,j} \in \mathcal{I}_{\mathbf{n}}}{\mathbf{i}\neq \mathbf{j}}}\Theta_{k\ell}^{(L)}(\lVert\mathbf{i}-\mathbf{j} \lVert_{2})=\theta_{k\ell}^{(L)},\,\,\,\lim_{L\rightarrow +\infty}\left(\theta_{k\ell}^{(L)}\right)=\theta_{k\ell}     \,\,\,\textrm{ and }\lim_{L\rightarrow +\infty}\left(\eta_{k\ell}^{(L)}\right)=\eta_{k\ell},
\end{equation}
where $\theta_{k\ell}^{(L)}$,  $\theta_{k\ell}$ and $\eta_{k\ell}$ are real numbers.
\end{itemize}
\end{assumption}
\begin{remark}  
\noindent Assumption \ref{ass0} is very usual in regression analysis where a stronger condition, given by equality of variances of the error terms,  is   often assumed. The $\alpha$-mixing condition in Assumption \ref{ass1} is classical in the literature on spatially dependent data (e.g.  \cite{bouka18}, \cite{bouka23}, \cite{dabo16},\cite{tran90}). Assumption \ref{ass2} and \ref{ass3}   refer to a type of process which has already been considered in the literature, in particular in\cite{gneiting2010} and \cite{leonenko2017}  with regard to Gaussian processes.  Two examples of classes of spatial correlation functions that satisfy the conditions $\Psi_{k,k} (0)=1$ and  $\sum\limits_{t=1}^{+\infty}t^{d-1}\Psi_{k,\ell}(t)<+\infty$  are:
\begin{itemize}
\item[$\bullet$] The powered exponential class:
$$\Psi_{k\ell}(t)=\exp\left[-\left(\dfrac{t}{b_{k\ell}}\right)^{b_0}\right],$$
where $b_0>0$ and $b_{k\ell}>0$.
\item[$\bullet$] The Mat\'ern class: 
$$\Psi_{k\ell}(t)=\frac{2^{1-c}}{\mathbf{G}(c)}a_{k\ell}^ct^{c}K_{c}(a_{k\ell}t),$$
where $c>0$, $a_{k\ell}>0$, $\mathbf{G}$ is the usual Gamma function  and $K_c$ is a modified Bessel function of the second kind (see \cite{gneiting2010}).
\end{itemize}
 Note that the sequence $$\bigg(n^{-d} \sum\limits_{\stackrel{\mathbf{i,j} \in \mathcal{I}_{\mathbf{n}}}{\mathbf{i}\neq \mathbf{j}}}\Psi_{k\ell}(\lVert\mathbf{i}-\mathbf{j} \lVert_{2})\bigg)_{n\in\mathbb{N}^\ast}$$  is bounded since
\[
\dfrac{1}{n^{d}} \sum\limits_{\stackrel{\mathbf{i,j} \in \mathcal{I}_{\mathbf{n}}}{\mathbf{i}\neq \mathbf{j}}}\Psi_{k\ell}(\lVert\mathbf{i}-\mathbf{j} \lVert_{2})\leqslant \sum_{t=1}^{+\infty}t^{d-1}\Psi_{k,\ell}(t)<+\infty.
\] 
So, it is not too strong to assume the first condition in \eqref{cvcor} and \eqref{cvcor2}. Furthermore,
by using the dominated convergence theorem, it is easily seen that the covariance operator of the truncated variables 
\[
V_1^L=\mathbb{E}\bigg(\left(Z_\mathbf{i}^{(L)}-\mathbb{E} (Z_\mathbf{i}^{(L)})\right)\otimes \left(Z_\mathbf{i}^{(L)}-\mathbb{E} (Z_\mathbf{i}^{(L)})\right)\bigg)
\]
converges to $V_1$ as $L\rightarrow +\infty$. So, it is natural to assume that elements related to $V_1^L$ converge to their analogues related to $V_1$ as $L\rightarrow +\infty$. It is the case for the eigenvalues and, under some conditions, for the eigenvectors (see \cite{ferre2003}). This justifies the two last conditions in \eqref{cvcor} of Assumption \ref{ass2}. The two last of conditions \eqref{cvcor2} are justified analogously from the convergence of the covariance operator of $Z_\mathbf{i}^{(L)}\otimes Z_\mathbf{i}^{(L)}$ to that of $Z\otimes Z$ as $L\rightarrow +\infty$.
\end{remark}

\subsection{Results}
\noindent The following results ensure consistency of the used estimators and give that of the proposed variable selection method. Considering the usual   norm of operators given by  $\parallel T\parallel_{\infty} =\sup\limits_{ x \in E }\dfrac{\parallel Tx \parallel_{F} }{\parallel x\parallel_{E}}$, where $T$ belongs to the space $\mathscr{L}(E,F)$ of operators from $E$ to $F$ , we first have:

\bigskip

\begin{proposition}\label{convop}
Under assumptions $1$ to $3$,  if    $\varphi(t)=O(t^{-\theta})$ with   $\theta> 2d$, then we have:		
	\begin{equation}\label{resv1}
		\mathbb{E}\Big(\parallel\widehat{V}_{1}^{(\mathbf{n})}- V_{1}\parallel_{\infty}^{2}\Big)=O(n^{-d}\log(n)) ,
	\end{equation}
	\begin{equation}\label{resv12}
		\mathbb{E}\Big(\parallel\widehat{V}_{12}^{(\mathbf{n})}- V_{12}\parallel_{\infty}^{2}\Big)=O(n^{-d}\log(n)). 
	\end{equation}	
\end{proposition}

\bigskip

\noindent An immediate  consequence of this proposition, obtained by using Markov inequality, is the almost complete convergences  of  $\widehat{V}_{1}^{(\mathbf{n})}$ and $\widehat{V}_{12}^{(\mathbf{n})}$  to $V_1$ and $V_{12}$ respectively, as $n\rightarrow +\infty$. This obviously  leads to the following result which ensures the consistency of the estimator of the criterion introduced above:

\bigskip

\begin{corollary}
Under assumptions $1$ to $3$,  if    $\varphi(t)=O(t^{-\theta})$ with   $\theta> 2d$, then  for all $K \subset I$,  the sequence  $\widehat{\xi}_{K}^{\mathbf{(n)}}$ almost surely converges to $\xi_K$ as $n\rightarrow +\infty$.
\end{corollary}

\bigskip

\noindent The following theorem gives the consistency of the introduced method of variable selection via a convergence of $\widehat{I}_1^{(\mathbf{n})}$ to $I_1$:

\bigskip

 \begin{theorem} \label{convest} 
Under assumptions $1$ to $4$,  if    $\varphi(t)=O(t^{-\theta})$ with   $\theta> 4d$, then:
\begin{equation}\label{coni1}	 
\lim\limits_{\ n \rightarrow +\infty }P\big(\widehat{I}_1^{(\mathbf{n})} = I_1 \big) = 1.
\end{equation}
\end{theorem}	
\section{Simulations}\label{sim}
This section presents the results of simulations made in order to assess the finite-sample behaviour of the proposed variable selection method and to compare it with existing ones.  We applied these methods to  simulated spatial data in $\mathbb{Z}^2$.  
\subsection{The simulated model}\label{simul}
\noindent Dealing with the case where  $p=6$, $q=1$ and $d=2$, we generated, as in \cite{dabo16},  each vector $X_{(i,j)}=\left(X_{(i,j)}^{(1)},X_{(i,j)}^{(2)},X_{(i,j)}^{(3)},X_{(i,j)}^{(4)},X_{(i,j)}^{(5)},X_{(i,j)}^{(6)}\right)^T$, $(i,j)\in\{1,\cdots,n\}^2$, such that
\[
X_{(i,j)}^{(k)} = D_{(i,j)} \times \dfrac{1}{\sqrt{500}} \sum\limits_{\ell=1}^{1000}\cos(w(1,\ell)\times i + w(2,\ell)\times j + q(\ell)\times t_{k} + r(\ell)),
\]
where $w(i,\ell)$ and  $q(u)$, $( i,\ell)\in\{1,2\}\times\{1,\cdots,1000\}$,  are  i.i.d from $N(0,0.25)$ and  independent to $r(\ell)$, $ \ell=1,\cdots,1000$, which are i.i.d  from uniform distribution on $[-\pi,\pi]$, $t_k=1+1.5(k-1)$, $k=1,\cdots,6$, and  
\[
 D_{(i,j)} = \dfrac{1}{n^2}\sum\limits_{(m,l)\in\{1,\cdots,n\}^2} \exp\Big(-\dfrac{\sqrt{(i-m)^2+(j-l)^2}}{a} \Big),
\]
with $a=2, 5, 10, 25 $.  This last term controls the spatial dependency which decreases as $a$ increases. The responses $Y_{(i,j)}$ were  then generated  according to the model
\[
Y_{(i,j)}= B X_{(i,j)} + \varepsilon_{(i,j)}
\]
where  $B=(3,5,4,6,0,0)$ and $\{\varepsilon_{(i,j)}\}_{(i,j)\in\{1,\cdots,n\}^2}$ is a centered  Gaussian random field  with spatial covariance function defined by 
\[
C(\mathbf{h})=\kappa^2\exp\bigg(\frac{\Vert\mathbf{h}\Vert_2^2}{9}\bigg),
\]
with $\kappa^2=1, 4, 9$.
\subsection{Choosing optimal tuning parameters}
\noindent Our procedure  depends on two tuning parameters  $\gamma$ and $\beta$ which  may have inflence on its performance. We chosed these parameters from the data by using a     cross validation procedure in order to minimize mean squared error of  prediction. More precisely, by using the lexico-graphic order, the sample is rewritten as $\{(X_{\mathbf{i}_i},Y_{\mathbf{i}_i})\}_{1\le i\le n^2}$. For each $\ell$ in $\{1,\cdots,n^2\}$, after removing the $\ell$-th observation   $(X_{\mathbf{i}_\ell},Y_{\mathbf{i}_\ell})$ from the sample,  we apply our method for selecting variable on the remaining  sample $\{(X_{\mathbf{i}_i},Y_{\mathbf{i}_i})\}_{ i\in\{1,\cdots, n^2\}-\{\ell\}}$ with a given value for $(\gamma,\beta)$ in $]0,1/2[^2$. This leads to   an estimate     $\widehat{I}_1^{-\ell}$ of $I_1$. Then, we consider the linear predictor of the $\ell$-th response from the explanatory variables that have been selected  as
\[
\widehat{Y}_{\gamma,\beta}^{(\ell)}=X_{\mathbf{i}_\ell}^TA_{\widehat{I}_1^{-\ell}}^T\left(A_{\widehat{I}_1^{-\ell}}X_{\mathbf{i}_\ell}X_{\mathbf{i}_\ell}^TA_{\widehat{I}_1^{-\ell}}^T\right)^{-1}A_{\widehat{I}_1^{-\ell}}X_{\mathbf{i}_\ell}Y_{\mathbf{i}_\ell},
\]
and   the  cross-validation index
\[
CV(\gamma,\beta)=\frac{1}{n^2}\sum_{\ell=1}^{n^2}\big(Y_{\mathbf{i}_\ell}-\widehat{Y}_{\gamma,\beta}^{(\ell)}\big)^2.
\]
Then, the optimal value   $(\gamma_{opt},\beta_{opt})$ of    $(\gamma,\beta)$  is  obtain by minimizing $CV(\gamma,\beta)$ over $]0,1/2[^2$, that is:
\begin{equation}\label{opti}
(\gamma_{opt},\beta_{opt})=\underset{(\gamma,\beta)\in ]0,1/2[^2}{\mathrm{argmin}}CV(\gamma,\beta).
\end{equation}

\subsection{Simulation strategy and results}
\noindent We simulated   500  independent replications  of  samples generated as indicated above, over which three  criteria  were  calculated  in order to assess performance   of the methods: (i) the average of mean squared error of prediction  after variable selection (MSE); (ii) the proportion of equality  (PE)  of the set    of  selected   variables to  the true set of relevant variables, that is   $\widehat{I}_1=I_1$;  (iii) the average of number of selected variables (NV). For computing mean squared errors of prediction   two independent data sets are generated for each replication: training data and test data, each with sample size $n^2=144$, $576$, $900$.  The training data is used for computing optimal values of the tuning parameters of our method, as indicated indicated above, and for selecting variables for all the methods considered for comparison.   The test data is used for computing MSE.  We compared  our method (OM)  to three methods of \cite{wang2009} consisting in minimizing penalized least squares, each being  related to a particular penalty function: SQUAD penalty, Hard thresholding penalty (Hard)  and $L^1$ penalty (LASSO). Our method  was performed by using the penalty functions $f(x)=\ln(x)^{-0.1}$ and $g(x)=\ln(x)^{0.1}$. Tables  \ref{tab:grav1} to  \ref{tab:grav3} report the obtained results. It can be seen that our method outperforms the three others since the values of   PE obtained for it are  greater than that obtained for SQUAD, Hard and LASSO. It can also be  observed that the weaker the spatial dependence the better the performance of this method. Indeed, PE has higher values for $a=25$ than for $a=5$ and $a=10$. Concerning  MSE, all the methods seem to be equivalent since the differences between the obtained values are very slight.

\begin{table}[htbp]
\centering 
	\begin{tabular}{cccccccccccc}
		\hline\hline
                           & & & &        & &  $n^2=144$ & &  & & &\\
\cline{2-12}
                          & &$a=5$& & & &$a=10$& &  & &$a=25$&\\
\cline{2-4}\cline{6-8}\cline{10-12}
		                       &MSE&NV&PE& &MSE&NV&PE&  &MSE&NV&PE \\
		$\kappa^2 = 1$ \\
		OM    				&0.125&3.60&0.386   & &0.109&4.300&0.402&    &0.120&4.600&0.344\\
		SCAD            &	0.074&5.300&0.068 &   &0.090&2.690&0.108 &    &0.109&5.480&0.056  \\
		Hard				&0.074&5.300&0.074 &   &0.090&5.410&0.052 &   & 0.109&5.420&0.088 \\
		LASSO           &	0.076&4.780&0.366   &  &0.092&5.110&0.158 &     &0.111&5.270&0.098  \\
		$\kappa^2 = 4$ \\
		OM           	&   0.171&3.730&0.184 &  &0.161&4.530&0.236& &	0.167&4.790&0.242\\
		SCAD            &  0.134&5.310&0.012  &  &0.145&5.380&0.024 &  &	0.159&5.450&0.016  \\
		Hard			&  0.134&5.360&0.022 &   &0.146&5.350&0.068 & &	0.159&5.380&0.054\\
		LASSO           &	 0.180&4.880&0.166  &  &0.172&5.260&0.158 &   &0.113&5.290&0.100 \\
		$\kappa^2 = 9$ \\
		OM    			 &0.226&3.780&0.244  & & 0.218&4.690&0.190  & &0.218&4.200&0.262\\
		SCAD               &0.195&5.350&0.020 &  &0.205&5.270&0.014&   &0.211&5.330&0.022 \\
		Hard			 &0.195&5.530&0.018 &  &0.205&5.340&0.030 &  &0.211&5.330&0.040  \\
		LASSO           &0.294&4.670&0.218& &0.209&5.070&0.120 &   &0.214&5.230&0.080 \\
\hline\hline
	\end{tabular}
	\caption{MSE, PE and NV  for different methods methods over $500$ replications, with $n=12$.}
	\label{tab:grav1} 
\end{table}

\begin{table}[htbp]
	\centering 
	\begin{tabular}{cccccccccccc}
		\hline\hline
                           & & & &        & &  $n^2=576$ & &  & & &\\
\cline{2-12}
                          & &$a=5$& & & &$a=10$& &  & &$a=25$&\\
\cline{2-4}\cline{6-8}\cline{10-12}
		                       &MSE&NV&PE& &MSE&NV&PE&  &MSE&NV&PE \\
		$\kappa^2 = 1$ \\
		OM   			& 	 0.079&1.190&0.000   &  &0.066&4.110&0.730&   & 0.046&4.080&0.684\\
		SCAD            &	 0.038&5.210&0.046&   &0.053&5.410&0.280&   & 0.046&5.460&0.098  \\
		Hard			&	 0.038&5.230&0.042&    &0.043&5.520&0.312&    &0.046&5.530&0.134\\
		LASSO           &	 0.039&3.900&0.394&   &0.041&4.640&0.512&    &0.046&5.110&0.184 \\		
		$\kappa^2 = 4$ \\		
		OM            	&    0.103&1.300&0.000&     &0.094&3.620&0.420&	 &0.081&4.340&0.602\\
		SCAD            &   0.076&5.210&0.008&     &0.077&5.340&0.052 &	 &0.080&5.560&0.056  \\
		Hard			&   0.076&5.320&0.014&    &0.077&5.370&0.088 &	&0.081&5.430&0.148\\
		LASSO           &	 0.077&3.810&0.290&   &0.077&4.690&0.460&   &0.080&5.100&0.188  \\
		$\kappa^2 = 9$ \\
		OM    			&0.135&1.380&0.012&   &0.131&3.690&0.382&  &0.120&4.110&0.312\\
		SCAD            &0.115&5.310&0.004 &   &0.116&5.270&0.028 &  &0.117&5.580&0.004  \\
		Hard			&0.115&5.460&0.008&   &0.116&5.230&0.076&   &0.117&5.230&0.096  \\
		LASSO           &0.116&3.770&0.210 & &0.116&4.750&0.348 &  &0.117&5.140&0.140 \\
\hline\hline
	\end{tabular}
	\caption{MSE, PE and NV  for different methods methods over $500$ replications, with $n=24$.}
	\label{tab:grav2} 
\end{table}
\begin{table}[htbp]
	\centering 
	\begin{tabular}{cccccccccccc}
		\hline\hline
                           & & & &        & &  $n^2=900$ & &  & & &\\
\cline{2-12}
                          & &$a=5$& & & &$a=10$& &  & &$a=25$&\\
\cline{2-4}\cline{6-8}\cline{10-12}
		                       &MSE&NV&PE& &MSE&NV&PE&  &MSE&NV&PE \\
		$\kappa^2 = 1$ \\
		OM    			& 	 0.054&1.000&0.000  &      &0.041&3.580&0.624&   & 0.039&4.020&0.720\\
		SCAD            &	 0.031&5.040&0.096 &   &0.032&5.380&0.184&    & 0.035&5.33&0.140  \\
		Hard			&	 0.031&5.090&0.060&    &0.032&5.540&0.104&    & 0.035&5.520&0.116\\
		LASSO           &	 0.032&3.040&0.204&   &0.032&4.350&0.680&    &0.035&4.99&0.204 \\
		$\kappa^2 = 4$ \\		
		OM            	&    0.078&1.000&0.000&        &0.091&3.710&0.442&	&0.073&4.310&0.630\\
		SCAD            &   0.063&5.320&0.024&   &0.075&5.210&0.074 &	&0.072&5.670&0.074 \\
		Hard			&   0.063&5.290&0.032&    &0.075&5.270&0.136 &	&0.072&5.550&0.164\\
		LASSO           &	 0.064&3.110&0.228&  &0.075&4.530&0.484&  &0.071&5.110&0.178 \\
		$\kappa^2 = 9$ \\		
		OM    			&0.125&1.420&0.012&   &0.129&3.820&0.370 & &0.131&4.090&0.340\\
		SCAD            &0.176&5.410&0.030&  &0.104&5.010&0.042&   &0.126&5.280&0.012  \\
		Hard			&0.177&5.510&0.022 &  &0.104&5.130&0.090&    &0.125&5.250&0.104  \\
		LASSO           &0.107&3.71&0.222&  &0.104&4.650&0.336 &   &0.126&5.050&0.206 \\
\hline\hline
	\end{tabular}
	\caption{MSE, PE and NV  for different methods methods over $500$ replications, with $n=30$.}
	\label{tab:grav3} 
\end{table}

\section{Proofs}\label{preuves}
\subsection{Proof of Proposition \ref{convop}}
\noindent Proof of   \eqref{resv1}: 
It is  proved in  \cite{bouka18} that 
\begin{equation}\label{bouka}
\mathbb{E}\big(\parallel n^{-d} \sum\limits_{\mathbf{i} \in \mathcal{I}_{\mathbf{n}}} X_{\mathbf{i}}\otimes X_{\mathbf{i}} - 	\mathbb{E}[X\otimes X] \parallel^{2}_{\mathcal{H}}\big)  = O(n^{-d}\log(n)),
\end{equation} 
and since
 $$	\mathbb{E}\big(\parallel \widehat{V}_{1}^{(\mathbf{n})}-V_{1} \parallel_{\mathcal{H}}^{2} \big) \leqslant 2	\mathbb{E}\big(\parallel n^{-d} \sum\limits_{\mathbf{i} \in \mathcal{I}_{\mathbf{n}}} X_{\mathbf{i}}\otimes X_{\mathbf{i}} - 	\mathbb{E}[X\otimes X] \parallel^{2}_{\mathcal{H}}\big) +  2	\mathbb{E}\big(\parallel\overline{X}^{(\mathbf{n})}\otimes \overline{X}^{(\mathbf{n}}\parallel^{2}_{\mathcal{H}}\big),$$ 
it remains to show that  
\begin{equation}\label{reste}
\mathbb{E}\big(\parallel\overline{X}^{(\mathbf{n})}\otimes \overline{X}^{(\mathbf{n}}\parallel^{2}_{\mathcal{H}}\big)= O(n^{-d}\log(n)).
\end{equation} 
Letting   $\{e_k\}_{1\leq k\leq p}$ be the canonical basis of $\mathbb{R}^p$,  putting $X_{\mathbf{i}}^{(k)}=\big<X_{\mathbf{i}},e_k\big>_{ \mathbb{R}^p}$ where  $\big<\cdot,\cdot\big>_{ \mathbb{R}^p}$ is the usual Euclidean inner product  of $\mathbb{R}^p$ and using the Cauchy-Schwartz inequality, we get:
\begin{eqnarray}\label{esp}
	\mathbb{E}\big(\parallel\overline{X}^{(\mathbf{n})}\otimes \overline{X}^{(\mathbf{n})}\parallel_{\mathcal{H}}^{2}\big)
	& = & \sum\limits_{j=1}^{p}\sum\limits_{k=1}^{p}	\mathbb{E}\Big(\big<(\overline{X}^{(\mathbf{n})}\otimes \overline{X}^{(\mathbf{n})})e_k,e_j\big>^{2}_{ \mathbb{R}^p}\Big)\nonumber\\
& = & \sum\limits_{j=1}^{p}\sum\limits_{k=1}^{p}	\mathbb{E}\Big(\big<\overline{X}^{(\mathbf{n})},e_k \big>^{2}_{ \mathbb{R}^p}\big< \overline{X}^{(\mathbf{n})},e_j\big>^{2}_{ \mathbb{R}^p}\Big)\nonumber\\
	 &\leq& p^{2} \mathbb{E}\big[ \parallel \overline{X}^{(\mathbf{n})} \parallel^{4}_{ \mathbb{R}^p} \big]\nonumber\\
	&=&n^{-4d}p^{2}	\mathbb{E}\bigg(\Big(\sum\limits_{k=1}^{p}\,\,\sum\limits_{(\mathbf{i}_{0},\mathbf{i}_{1}) \in \mathcal{I}_{\mathbf{n}}^2}   X_{\mathbf{i}_{0}}^{(k)}X_{\mathbf{i}_{1}}^{(k)}\Big)^{2}\bigg)\nonumber\\
	&=& n^{-4d}p^{2}\sum\limits_{k=1}^{p}\sum\limits_{i=1}^{p}\,\,\,\sum\limits_{(\mathbf{i}_{0},\mathbf{i}_{1},\mathbf{i}_{2},\mathbf{i}_{3}) \in \mathcal{I}_{\mathbf{n}}^4} 	\mathbb{E}\big(X_{\mathbf{i}_{0}}^{(k)}X_{\mathbf{i}_{1}}^{(k)}X_{\mathbf{i}_{2}}^{(i)}X_{\mathbf{i}_{3}}^{(i)} \big).
\end{eqnarray}
Let us  denote by  $\mathcal{J}_{\mathbf{n}}^{(\ell)}$, $\ell= 1,2,3,4$, the subsets of  $\mathcal{I}_{\mathbf{n}}^4$ defined as follows:
\begin{itemize}
\item[$\bullet$] $\mathcal{J}_{\mathbf{n}}^{(1)}=\big\{(\mathbf{i}_{0}, \mathbf{i}_{1},\mathbf{i}_{2},\mathbf{i}_{3})\in \mathcal{I}_{\mathbf{n}}^4\,\,/\,\,{\mathbf{i}_{0}}={\mathbf{i}_{1}}={\mathbf{i}_{2}}={\mathbf{i}_{3}}\big\}$;
\item[$\bullet$] $\mathcal{J}_{\mathbf{n}}^{(2)}$ consists of the terms $(\mathbf{i}_{0}, \mathbf{i}_{1},\mathbf{i}_{2},\mathbf{i}_{3})\in \mathcal{I}_{\mathbf{n}}^4$ for which  one of the following properties holds: ${\mathbf{i}_{0}}={\mathbf{i}_{2}}\neq{\mathbf{i}_{1}}={\mathbf{i}_{3}}$, ${\mathbf{i}_{0}}={\mathbf{i}_{1}}\neq{\mathbf{i}_{2}}={\mathbf{i}_{3}}$, ${\mathbf{i}_{0}}={\mathbf{i}_{3}}\neq{\mathbf{i}_{1}}={\mathbf{i}_{2}}$, ${\mathbf{i}_{0}}\neq{\mathbf{i}_{1}}={\mathbf{i}_{2}}={\mathbf{i}_{3}}$, ${\mathbf{i}_{1}}\neq{\mathbf{i}_{0}}={\mathbf{i}_{2}}={\mathbf{i}_{3}}$, ${\mathbf{i}_{2}}\neq{\mathbf{i}_{0}}={\mathbf{i}_{1}}={\mathbf{i}_{3}}$ and ${\mathbf{i}_{3}}\neq{\mathbf{i}_{0}}={\mathbf{i}_{1}}={\mathbf{i}_{2}}$;
\item[$\bullet$] $\mathcal{J}_{\mathbf{n}}^{(3)}$ consists of the terms $(\mathbf{i}_{0}, \mathbf{i}_{1},\mathbf{i}_{2},\mathbf{i}_{3})\in \mathcal{I}_{\mathbf{n}}^4$ for which  one of the following properties holds: ${\mathbf{i}_{0}}={\mathbf{i}_{2}}\neq{\mathbf{i}_{1}}\neq{\mathbf{i}_{3}}$, ${\mathbf{i}_{0}}={\mathbf{i}_{1}}\neq{\mathbf{i}_{2}}\neq{\mathbf{i}_{3}}$, ${\mathbf{i}_{0}}={\mathbf{i}_{3}}\neq{\mathbf{i}_{1}}\neq{\mathbf{i}_{2}}$, ${\mathbf{i}_{1}}={\mathbf{i}_{2}}\neq{\mathbf{i}_{0}}\neq{\mathbf{i}_{3}}$, ${\mathbf{i}_{1}}={\mathbf{i}_{3}}\neq{\mathbf{i}_{0}}\neq{\mathbf{i}_{2}}$ and ${\mathbf{i}_{2}}={\mathbf{i}_{3}}\neq{\mathbf{i}_{0}}\neq{\mathbf{i}_{1}}$;
\item[$\bullet$] $\mathcal{J}_{\mathbf{n}}^{(4)}=\big\{(\mathbf{i}_{0}, \mathbf{i}_{1},\mathbf{i}_{2},\mathbf{i}_{3})\in \mathcal{I}_{\mathbf{n}}^4\,\,/\,\,\mathbf{i}_{0}\neq  \mathbf{i}_{1}, \mathbf{i}_{0}\neq  \mathbf{i}_{2}, \mathbf{i}_{0}\neq  \mathbf{i}_{3}, \mathbf{i}_{1}\neq  \mathbf{i}_{2}, \mathbf{i}_{1}\neq  \mathbf{i}_{3}, \mathbf{i}_{2}\neq  \mathbf{i}_{3}\big\}$.
\end{itemize}
Then, clearly, 
\begin{equation}\label{decomp}
\sum\limits_{(\mathbf{i}_{0},\mathbf{i}_{1},\mathbf{i}_{2},\mathbf{i}_{3}) \in \mathcal{I}_{\mathbf{n}}^4} 	\mathbb{E}\big(X_{\mathbf{i}_{0}}^{(k)}X_{\mathbf{i}_{1}}^{(k)}X_{\mathbf{i}_{2}}^{(i)}X_{\mathbf{i}_{3}}^{(i)} \big)=A_{\mathbf{n}}^{(1)}+A_{\mathbf{n}}^{(2)}+A_{\mathbf{n}}^{(3)}+A_{\mathbf{n}}^{(4)},
\end{equation}
where, for $\ell\in\{1,2,3,4\}$,
\[
A_{\mathbf{n}}^{(\ell)}=\sum\limits_{(\mathbf{i}_{0},\mathbf{i}_{1},\mathbf{i}_{2},\mathbf{i}_{3}) \in \mathcal{J}_{\mathbf{n}}^{(\ell)}} 	\mathbb{E}\big(X_{\mathbf{i}_{0}}^{(k)}X_{\mathbf{i}_{1}}^{(k)}X_{\mathbf{i}_{2}}^{(i)}X_{\mathbf{i}_{3}}^{(i)} \big).
\]
By using twice the  H$\ddot{\textrm{o}}$lder's inequality, together with   the strict stationarity property (see Assumption \ref{ass1}), it follows:
\begin{eqnarray*}
\mathbb{E}\big(\big\vert X_{\mathbf{i}_{0}}^{(k)}X_{\mathbf{i}_{1}}^{(k)}X_{\mathbf{i}_{2}}^{(i)}X_{\mathbf{i}_{3}}^{(i)}\big\vert \big)
&\leq &\mathbb{E}\big(\big( X_{\mathbf{i}_{0}}^{(i)}\big)^{4}\big)^{1/4}\mathbb{E}\big(\big(X_{\mathbf{i}_{1}}^{(k)}\big)^{4}\big)^{1/4}	\mathbb{E}\big(\big( X_{\mathbf{i}_{2}}^{(i)}\big)^{4}\big)^{1/4}\mathbb{E}\big(\big(X_{\mathbf{i}_{3}}^{(k)}\big)^{4}\big)^{1/4}=\mathbb{E}\big(X^4\big);
\end{eqnarray*}
hence
\[
\big\vert A_{\mathbf{n}}^{(1)}\big\vert\leq \mathbb{E}\big(X^4\big) \sum\limits_{\mathbf{i}_{0}\in \mathcal{I}_{\mathbf{n}}}^{}1=\mathbb{E}\big(X^4\big)n^d,
\]
\[
\big\vert A_{\mathbf{n}}^{(2)}\big\vert\leq 7\,\, \mathbb{E}\big(X^4\big)\sum\limits_{\stackrel{(\mathbf{i}_{0},\mathbf{i}_{1}) \in\mathcal{I}_{\mathbf{n}}^2}{ \mathbf{i}_{0}\neq\mathbf{i}_{1}} }  1\leq  7\,\, \mathbb{E}\big(X^4\big)\sum_{(\mathbf{i}_{0},\mathbf{i}_{1}) \in\mathcal{I}_{\mathbf{n}}^2}  1=7\mathbb{E}\big(X^4\big)n^{2d}
\]
and
\[
\big\vert A_{\mathbf{n}}^{(3)}\big\vert\leq 6\,\, \mathbb{E}\big(X^4\big)\sum\limits_{\stackrel{(\mathbf{i}_{0},\mathbf{i}_{1},\mathbf{i}_{2}) \in\mathcal{I}_{\mathbf{n}}^3}{ \mathbf{i}_{0}\neq\mathbf{i}_{1}\neq\mathbf{i}_{2}} }  1\leq  6\,\, \mathbb{E}\big(X^4\big)\sum_{(\mathbf{i}_{0},\mathbf{i}_{1},\mathbf{i}_{2}) \in\mathcal{I}_{\mathbf{n}}^3}  1=6\mathbb{E}\big(X^4\big)n^{3d},
\]
therefore
\begin{equation}\label{grando}
A_{\mathbf{n}}^{(1)}= O(n^{d}),\,\,\,A_{\mathbf{n}}^{(2)}= O(n^{2d})\,\,\,\textrm{ and }\,\,\,A_{\mathbf{n}}^{(3)}= O(n^{3d}).
\end{equation}
For the last term, by a similar reasoning than in \cite{gao2008} (see Equation (3.20)), we get
\begin{eqnarray*}
\sum\limits_{(\mathbf{i}_{0},\mathbf{i}_{1},\mathbf{i}_{2},\mathbf{i}_{3}) \in \mathcal{J}_{\mathbf{n}}^{(4)}} \Big\vert	\mathbb{E}\big(X_{\mathbf{i}_{0}}^{(k)}X_{\mathbf{i}_{1}}^{(k)}X_{\mathbf{i}_{2}}^{(i)}X_{\mathbf{i}_{3}}^{(i)} \big)\Big\vert
&\leqslant & C_{1}n^{2d}C_{n}^{d} + C_{2}n^{2d}\sum\limits_{t=1}^{+\infty}t^{d-1}\varphi(t)^{1/2}\\
&\leqslant&  C_{1}n^{2d}C_{n}^{d} + C_{2}n^{2d}\sum\limits_{t=1}^{+\infty}t^{d-1-\theta/2},
\end{eqnarray*}
where $C_1$ and $C_2$ are positive constants. Then, $A_{\mathbf{n}}^{(4)}= O(n^{2d})$ and from \eqref{grando} and \eqref{decomp} it follows that
\[
\sum\limits_{(\mathbf{i}_{0},\mathbf{i}_{1},\mathbf{i}_{2},\mathbf{i}_{3}) \in \mathcal{I}_{\mathbf{n}}^4} 	\mathbb{E}\big(X_{\mathbf{i}_{0}}^{(k)}X_{\mathbf{i}_{1}}^{(k)}X_{\mathbf{i}_{2}}^{(i)}X_{\mathbf{i}_{3}}^{(i)} \big)=O(n^{3d}).
\]
An use of  \eqref{esp} yields   $\mathbb{E}\big(\parallel\overline{X}^{(\mathbf{n})}\otimes \overline{X}^{(\mathbf{n}}\parallel^{2}_{\mathcal{H}}\big)=O(n^{-d})= O(n^{-d}\log(n))$, what allows to conclude that  $\mathbb{E}\big(\Vert\widehat{V}_{1}^{(\mathbf{n})}-V_{1} \Vert_{\mathcal{H}}^{2} \big) =O(n^{-d}\log(n))$. Then \eqref{resv1} is  obtained from $\Vert\cdot\Vert_\infty\leqslant \Vert\cdot\Vert_\mathcal{H}$.

\bigskip

\noindent Proof of   \eqref{resv12}: We have:
\begin{eqnarray}\label{v12}
	\mathbb{E}\Big(\parallel \widehat{V}_{12}-V_{12}\parallel_{\infty}^{2} \Big) & \leqslant & 2 \mathbb{E}\Big(  \Big\lVert n^{-d} \sum\limits_{\mathbf{i} \in \mathcal{I}_{\mathbf{n}}} Y _{\mathbf{i}}\otimes X_{\mathbf{i}} - \mathbb{ E}(Y \otimes X) \Big\lVert_{\infty}^{2} \Big)+ 2\mathbb{E}\Big(\Big\lVert\overline{Y }^{(\mathbf{n})}\otimes \overline{X}^{(\mathbf{n})}\Big\lVert_{\infty}^{2}\Big)\nonumber\\
	& := &D_{1} + D_{2}.
\end{eqnarray}
Using \eqref{model} and a property of tensor products (see \cite{dauxois94}), we get
\begin{eqnarray*}
  Y _{\mathbf{i}}\otimes X_{\mathbf{i}} - \mathbb{ E}(Y \otimes X) &=&(BX _{\mathbf{i}})\otimes X_{\mathbf{i}} - \mathbb{ E}((BX) \otimes X)+\varepsilon_{\mathbf{i}}\otimes X_{\mathbf{i}}-\mathbb{E}(\varepsilon \otimes  X)\\
 &=&\Big(X _{\mathbf{i}}\otimes X_{\mathbf{i}} - \mathbb{ E}(X \otimes X)\Big)B^\ast+\varepsilon_{\mathbf{i}}\otimes X_{\mathbf{i}}-\mathbb{E}(\varepsilon )\otimes  \mathbb{E}(X)\\
 &=&\Big(X _{\mathbf{i}}\otimes X_{\mathbf{i}} - \mathbb{ E}(X \otimes X)\Big)B^\ast+\varepsilon_{\mathbf{i}}\otimes X_{\mathbf{i}}
\end{eqnarray*}
and, therefore,
\begin{eqnarray*}
	D_{1} & \leqslant & 4\lVert B\lVert_{\infty}^{2}\mathbb{E}\Big(  \Big\lVert n^{-d} \sum\limits_{\mathbf{i} \in \mathcal{I}_{\mathbf{n}}} X_{\mathbf{i}}\otimes X_{\mathbf{i}} -  \mathbb{E}(X\otimes X) \Big\lVert_{\mathcal{H}}^{2} \Big) + 4\mathbb{E}\Big( \Big\lVert n^{-d} \sum\limits_{\mathbf{i} \in \mathcal{I}_{\mathbf{n}}} \varepsilon_{\mathbf{i}}\otimes X_{\mathbf{i}} \Big\lVert_{\mathcal{H}}^{2} \Big).
\end{eqnarray*}
Using  the properties of the tensor  product (see \cite{dauxois94}), in particular the   equality $\left\langle a\otimes b,c\otimes d\right\rangle_\mathcal{H}=\langle a,c\rangle\,\langle b,d\rangle$, we obtain
\begin{eqnarray*}
	\mathbb{E}\Big( \Big\lVert \sum\limits_{\mathbf{i} \in \mathcal{I}_{\mathbf{n}}} \varepsilon_{\mathbf{i}}\otimes X_{\mathbf{i}} \Big\lVert_{\mathcal{H}}^{2} \Big)&=&\sum\limits_{(\mathbf{i},\mathbf{j}) \in \mathcal{I}_{\mathbf{n}}^2}\mathbb{E}\Big( \Big\langle \varepsilon_{\mathbf{i}}\otimes X_{\mathbf{i}},\varepsilon_{\mathbf{j}}\otimes X_{\mathbf{j}}\Big\rangle_\mathcal{H}\Big)\\
&=&\sum\limits_{(\mathbf{i},\mathbf{j}) \in \mathcal{I}_{\mathbf{n}}^2}\mathbb{E}\Big( \Big\langle X_{\mathbf{i}}, X_{\mathbf{j}}\Big\rangle_{ \mathbb{R}^p} \Big\langle\varepsilon_{\mathbf{i}}, \varepsilon_{\mathbf{j}}\Big\rangle_{ \mathbb{R}^q}\Big)\\
&=&\sum\limits_{(\mathbf{i},\mathbf{j}) \in \mathcal{I}_{\mathbf{n}}^2}\mathbb{E}\Big( \Big\langle X_{\mathbf{i}}, X_{\mathbf{j}}\Big\rangle_{ \mathbb{R}^p} \Big\langle\varepsilon_{\mathbf{i}}, \varepsilon_{\mathbf{j}}\Big\rangle_{ \mathbb{R}^q}\Big)\\
&=&\sum\limits_{\stackrel{(\mathbf{i},\mathbf{j}) \in \mathcal{I}_{\mathbf{n}}^2}{\mathbf{i}\neq \mathbf{j}}}\mathbb{E}\Big( \Big\langle X_{\mathbf{i}}, X_{\mathbf{j}}\Big\rangle_{ \mathbb{R}^p} \Big\langle\varepsilon_{\mathbf{i}}, \varepsilon_{\mathbf{j}}\Big\rangle_{ \mathbb{R}^q}\Big)+ 
\sum\limits_{\mathbf{i} \in \mathcal{I}_{\mathbf{n}}}\mathbb{E}\Big( \left\Vert X_{\mathbf{i}}\right\Vert_{ \mathbb{R}^p} ^2 \left\Vert \varepsilon_{\mathbf{i}} \right\Vert _{ \mathbb{R}^q}^2\Big).
\end{eqnarray*}
Using this equality together with  \eqref{bouka}, the decomposition $X_{\mathbf{i}}=\sum_{k=1}^p\langle X_{\mathbf{i}},u_k\rangle_{ \mathbb{R}^p} u_k$, where $\{u_k\}_{1\leq k\leq p}$ is an orthonormal basis of $\mathbb{R}^p$ consisting of eigenvectors of $V_1$ such that $u_k$ is associated with the $k$-th largest eigenvalue $\lambda_k$, and Assumptions 1 and 3,  it follows
	\begin{eqnarray*}
	D_1& \leqslant & C_1n^{-d}\log(n) +   4n^{-2d}\sum\limits_{\ell=1}^{p}\sum\limits_{k=1}^{p}\sum\limits_{\stackrel{(\mathbf{i},\mathbf{j}) \in \mathcal{I}_{\mathbf{n}}^2}{\mathbf{i}\neq \mathbf{j}}}\mathbb{E}\Big(\Big<X_{\mathbf{i}},u_\ell\Big>_{ \mathbb{R}^p}\Big<X_{\mathbf{j}},u_k\Big>_{ \mathbb{R}^p} \Big)\mathbb{E}\Big(\Big<\varepsilon_{\mathbf{i}},\varepsilon_{\mathbf{j}}\Big>_{ \mathbb{R}^q} \Big) \\
	 & + & 4n^{-2d}\sum\limits_{\mathbf{i} \in \mathcal{I}_{\mathbf{n}}}\mathbb{E}\Big( \parallel X_{\mathbf{i}}\parallel^{2} _{ \mathbb{R}^p}\Big)\mathbb{E}\Big( \parallel \varepsilon_{\mathbf{i}}\parallel^{2} _{ \mathbb{R}^q}\Big)\\
	 & \leqslant & C_1n^{-d}\log(n) \\
&+ &4n^{-2d}\sum\limits_{\ell=1}^{p}\sum\limits_{k=1}^{p}\sum\limits_{\stackrel{(\mathbf{i},\mathbf{j}) \in \mathcal{I}_{\mathbf{n}}^2}{\mathbf{i}\neq \mathbf{j}}} \Big|\mathbb{E}\Big(\langle X_{\mathbf{i}},u_\ell\rangle_{ \mathbb{R}^p}\langle X_{\mathbf{j}},u_k\rangle_{ \mathbb{R}^p}\Big)\Big| \mathbb{E}\big(\parallel \varepsilon_{\mathbf{i}}\parallel^{2}_{ \mathbb{R}^q}\big)^{1/2} \mathbb{E}\big(\parallel \varepsilon_{\mathbf{j}}\parallel^{2}_{ \mathbb{R}^q}\big)^{1/2} \\
	  & + & C_2n^{-d}\\
	   & \leqslant & C_1n^{-d}\log(n) +n^{-d}\sum\limits_{k=1}^{p}\sum\limits_{\ell=1}^{p}\left\vert\rho_{k\ell}\right\vert\lambda_{k}^{1/2}\lambda_{\ell}^{1/2}\sum\limits_{t=1}^{+\infty}t^{d-1}\Psi_{k,\ell}(t)+C_2n^{-d},
\end{eqnarray*}
where $C_1>0$ and  $C_2>0$. This implies that $D_{1} = O(n^{-d}\log(n))$. On the other hand, putting $\overline{\varepsilon}^{(\mathbf{n})}=n^{-d}\sum\limits_{\mathbf{i} \in \mathcal{I}_{\mathbf{n}}}\varepsilon_{\mathbf{i}}$ and using  as above  \eqref{model} and properties of the tensor product,   we obtain
\begin{eqnarray}\label{d2}
	D_{2} & \leqslant &  4\lVert B\lVert_{\infty}^{2}\mathbb{E}\Big(\Big\lVert\overline{X}^{(\mathbf{n})}\otimes \overline{X}^{(\mathbf{n})}\Big\lVert_{\mathcal{H}}^{2}\Big) + 4\mathbb{E}\Big(\Big\lVert \overline{\varepsilon}^{(\mathbf{n})} \otimes \overline{X}^{(\mathbf{n})}\Big\lVert_{\mathcal{H}}^{2}\Big)\nonumber\\
&=&  4\lVert B\lVert_{\infty}^{2}\mathbb{E}\Big(\Big\lVert\overline{X}^{(\mathbf{n})}\otimes \overline{X}^{(\mathbf{n})}\Big\lVert_{\mathcal{H}}^{2}\Big) + 4\mathbb{E}\big(\parallel\overline{X}^{(\mathbf{n})} \parallel^{2}_{ \mathbb{R}^p}\big) \mathbb{E}\big(\parallel\overline{\mathbf{\varepsilon}}^{(\mathbf{n})} \parallel^{2}_{ \mathbb{R}^q}\big).
\end{eqnarray}
Since
\begin{equation}\label{const}
 \mathbb{E}\big(\parallel\overline{\mathbf{\varepsilon}}^{(\mathbf{n})} \parallel^{2}_{ \mathbb{R}^q}\big)=\frac{1}{n^{2d}}\sum _{(\mathbf{i},\mathbf{j}) \in \mathcal{I}_{\mathbf{n}}^2}\mathbb{E}\Big(\left\langle \varepsilon_{\mathbf{i}}, \varepsilon_{\mathbf{j}}\right\rangle_{ \mathbb{R}^q}\Big)\leq  \frac{1}{n^{2d}}\sum _{(\mathbf{i},\mathbf{j}) \in \mathcal{I}_{\mathbf{n}}^2}\mathbb{E}\Big(\Vert \varepsilon_{\mathbf{i}}\Vert^2_{ \mathbb{R}^q}\Big)^{1/2}\mathbb{E}\Big(\Vert \varepsilon_{\mathbf{j}}\Vert^2_{ \mathbb{R}^q}\Big)^{1/2}\leqslant\sigma^2_\varepsilon,
\end{equation}
where $\sigma^2_\varepsilon=\sup_{\mathbf{i}\in\mathbb{Z}^d}\mathbb{E}\left(\Vert \varepsilon_\mathbf{i}\Vert_{\mathbb{R}^q}^2\right)$, and  
\begin{eqnarray}\label{dern}
& &\mathbb{ E}\big(\parallel\overline{X}^{(\mathbf{n})} \parallel^{2}_{ \mathbb{R}^p}\big)\nonumber\\
&=&\mathbb{E}\Big(\Big\lVert\frac{1}{n^{d}}\sum_{k=1}^p\bigg(\sum _{\mathbf{i} \in \mathcal{I}_{\mathbf{n}}}\langle X_{\mathbf{i}},u_k\rangle_{ \mathbb{R}^p}\bigg)u_k\Big\lVert^{2}_{ \mathbb{R}^p}\Big)=\frac{1}{n^{2d}}\sum_{k=1}^p\mathbb{E}\Bigg(\bigg(\sum _{\mathbf{i} \in \mathcal{I}_{\mathbf{n}}}\langle X_{\mathbf{i}},u_k\rangle_{ \mathbb{R}^p}\bigg)^2\Bigg)\nonumber\\
&=&\frac{1}{n^{2d}}\sum_{k=1}^p\sum _{(\mathbf{i},\mathbf{j}) \in \mathcal{I}_{\mathbf{n}}^2}\mathbb{E}\Bigg(\langle X_{\mathbf{i}},u_k\rangle_{ \mathbb{R}^p} \langle X_{\mathbf{j}},u_k\rangle_{ \mathbb{R}^p}\Bigg)\nonumber\\
&=&n^{-2d}\sum\limits_{k=1}^{p}\sum\limits_{\stackrel{(\mathbf{i},\mathbf{j}) \in \mathcal{I}_{\mathbf{n}}^2}{\mathbf{i}\neq \mathbf{j}}}\mathbb{E}\Big(\Big<X_{\mathbf{i}},u_k\Big>_{ \mathbb{R}^p}\Big<X_{\mathbf{j}},u_k\Big> _{ \mathbb{R}^p}\Big)+n^{-2d}\sum\limits_{k=1}^{p}\sum\limits_{\mathbf{i} \in \mathcal{I}_{\mathbf{n}}}\mathbb{E}\Big( \Big<  X_{\mathbf{i}},u_k\Big>^{2}_{ \mathbb{R}^p} \Big)\nonumber\\
&\leqslant&n^{-2d}\sum\limits_{k=1}^{p}\sum\limits_{\stackrel{(\mathbf{i},\mathbf{j}) \in \mathcal{I}_{\mathbf{n}}^2}{\mathbf{i}\neq \mathbf{j}}} \Big|\mathbb{E}\Big(\langle X_{\mathbf{i}},u_k\rangle_{ \mathbb{R}^p}\langle X_{\mathbf{j}},u_k\rangle_{ \mathbb{R}^p}\Big)\Big| +n^{-2d}\sum\limits_{k=1}^{p}\sum\limits_{\mathbf{i} \in \mathcal{I}_{\mathbf{n}}}\mathbb{E}\Big( \parallel X_{\mathbf{i}}\parallel^{2}  _{ \mathbb{R}^p}\Big)\nonumber\\
&\leqslant&C_3 n^{-d}\sum\limits_{k=1}^{p}\lambda_{k}\sum\limits_{t=1}^{+\infty}t^{d-1}\Psi_{k,k}(t)+p\mathbb{E}(\Vert X\Vert^2_{ \mathbb{R}^p})n^{-d},
\end{eqnarray}
where $C_3>0$, it follows from \eqref{d2}, \eqref{reste}, \eqref{const} and \eqref{dern}  that  $D_2= O(n^{-d}\log(n))$,  and the required result  is deduced from \eqref{v12}.
\subsection{A central limit theorem}

\noindent This section is devoted to a central limit theorem which will be useful for proving the main results of the paper.
\begin{theorem}\label{clt}
Under Assumption 2, if $\varphi(t)=O(t^{-\theta})$ with $\theta>4d$  then, for all measurable and bounded function $h\,:\,\mathbb{R}^{p+q}\rightarrow\mathbb{R}$ such that there exists $\nu\in\mathbb{R}$ satisfying
\begin{equation}\label{condvar}
\lim_{n\rightarrow +\infty}\bigg\{\dfrac{1}{n^{d}}\sum_{    \stackrel{(\mathbf{i},\mathbf{j} )\in \mathcal{I}_{\mathbf{n}}^2}{  \mathbf{i}\neq \mathbf{j}  } }Cov\Big(h\left(Z_\mathbf{i}\right),h\left(Z_\mathbf{j}\right)\Big)\bigg\}=\nu,
\end{equation}
the random variable  $n^{-d/2}\sum\limits_{\mathbf{i} \in \mathcal{I}_{\mathbf{n}}}\Big(h\left(Z_\mathbf{i}\right)- \mathbb{E}(h(Z_\mathbf{i}))\Big) $ converges in distribution, as $n\rightarrow +\infty$, to a random variable   having a normal distribution $N(0,\sigma^2_h)$, where $\sigma^2_h=Var\left(h(Z)\right)+\nu$. 
\end{theorem}

\noindent\textit{Proof}. Clearly, $n^{-d/2}\sum\limits_{\mathbf{i} \in \mathcal{I}_{\mathbf{n}}}\Big(h\left(Z_\mathbf{i}\right)- \mathbb{E}(h(Z_\mathbf{i}))\Big) =n^{d/2}S_{\mathbf{n}} $, where $S_{\mathbf{n}} = \sum\limits_{\mathbf{i} \in \mathcal{I}_{\mathbf{n}}}\Delta_{\mathbf{i}}
$
with  $\Delta_{\mathbf{i}} =n^{-d}\Big( h\left(Z_\mathbf{i}\right)- \mathbb{E}(h(Z_\mathbf{i})) \Big)$.
Suppose, without loss of generality, that there exist  integers $(r,p_{1},p_{2})\in (\mathbb{N}^{*})^3$ such that  $n=r(p_{1} +p_{2})$ and $p_{2} < p_{1}$.  Let $p_{1} = \lfloor n^{1/3}\rfloor$ and $ p_{2}=\lfloor n^{1/4}\rfloor$, where $\lfloor\cdot\rfloor$ denotes the integer part. Then $p_{2}p_{1}^{-1} \leqslant Cn^{-1/12}$ and $p_{2} < p_{1}$ for $n$ large enough.  As in \cite{tran90}, the $\Delta_{\mathbf{i}}$s are set into large blocks and small blocks. For $\mathbf{j}=(j_{1},...,j_{d}) \in \{0,...,r-1\}^{d}$, 
let 
\begin{eqnarray*}
	U(1,\mathbf{n},\mathbf{j}) & = & \sum\limits_{\stackrel{i_{k}=j_{k}(p_{1}+p_{2})+1}{k=1,...,d}}^{j_{k}(p_{1}+p_{2})+p_{1}}\Delta_{\mathbf{i}},\\
	U(2,\mathbf{n},\mathbf{j}) & = & \sum\limits_{\stackrel{i_{k}=j_{k}(p_{1}+p_{2})+1}{k=1,...,d-1}}^{j_{k}(p_{1}+p_{2} )+p_{1}}\,\,\,\sum\limits_{{i_{d}=j_{d}(p_{1}+p_{2})+p_{1}+1}}^{(j_{d}+1)(p_{1}+p_{2})}\Delta_{\mathbf{i}},\\
	U(3,\mathbf{n},\mathbf{j}) & = & \sum\limits_{\stackrel{i_{k}=j_{k}(p_{1}+p_{2})+1}{k=1,...,d-2}}^{j_{k}(p_{1}+p_{2}) +p_{1}}\,\,\,\sum\limits_{{i_{d-1}=j_{d-1}(p_{1}+p_{2})+p_{1}+1}}^{(j_{d-1}+1)(p_{1}+p_{2})}\,\,\,\sum\limits_{i_{d}=j_{d}(p_{1}+p_{2})+1}^{j_{d}(p_{1} +p_{2})+p_{1}}\Delta_{\mathbf{i}},\\
	U(4,\mathbf{n},\mathbf{j}) & = & \sum\limits_{\stackrel{i_{k}=j_{k}(p_{1}+p_{2})+1}{k=1,...,d-2}}^{j_{k}(p_{1}+p_{2}) +p_{1}}\,\,\,\sum\limits_{{i_{d-1}=j_{d-1}(p_{1}+p_{2})+p_{1}+1}}^{(j_{d-1}+1)(p_{1}+p_{2})}\,\,\,\sum\limits_{i_{d}=j_{d}(p_{1}+p_{2})+p_{1}+1}^{(j_{d}+1)(p_{1}+p_{2})}\Delta_{\mathbf{i}},
\end{eqnarray*}
and so on. Finally
\begin{eqnarray*}
	U(2^{d-1},\mathbf{n},\mathbf{j}) & = & \sum\limits_{\stackrel{i_{k}=j_{k}(p_{1}+p_{2})+ p_{1}+1}{k=1,...,d-1}}^{(j_{k}+1)(p_{1}+p_{2})}\,\,\,\sum\limits_{{i_{d}=j_{d}(p_{1}+p_{2})+1}}^{j_{d}(p_{1}+p_{2})+p_{1}}\Delta_{\mathbf{i}},\\
	U(2^{d},\mathbf{n},\mathbf{j}) & = & \sum\limits_{\stackrel{i_{k}=j_{k}(p_{1}+p_{2})+1}{k=1,...,d}}^{(j_{k}+1)(p_{1}+p_{2})}\Delta_{\mathbf{i}}.
\end{eqnarray*}
The random variable $S_{\mathbf{n}} $ can be writen by using the preceding terms as 
$ S_{\mathbf{n}} = \sum\limits_{\ell=1}^{2^{d}} \mathscr{S}(\mathbf{n},\ell),$
where
$  \mathscr{S}(\mathbf{n},\ell) = \sum\limits_{\mathbf{j} \in \mathcal{J}}^{}U(\ell,\mathbf{n},\mathbf{j}),$
with   $\mathcal{J} = \{0,...,r-1 \}^{d}$. Then, in order to obtain asymptotic normality,   it is enough to show that 
\begin{equation}\label{cd1}
Q_{1}=\Bigg\vert \mathbb{E}\bigg[\exp\big(i\,n^{d/2}\mathscr{S}(\mathbf{n},1)\big)\bigg] - \prod\limits_{\mathbf{j} \in \mathcal{J}} \mathbb{E}\bigg[\exp\big(i\,n^{d/2}U(1,\mathbf{n},\mathbf{j}) \big) \bigg]\Bigg\vert \longrightarrow  0\,\,\,\textrm{ where }\,\,\,  i^{2}=-1  ,
\end{equation}
\begin{equation}\label{c2}
		 n^{d}\mathbb{E}\bigg[\Big( \sum\limits_{\ell=2}^{2^{d}}   \mathscr{S}(\mathbf{n},\ell)\Big)^{2}\bigg] \longrightarrow  0,
\end{equation}
\begin{equation}\label{c3}
	 n^{d}\sum\limits_{\mathbf{j} \in \mathcal{J}} \mathbb{E}\bigg[ U(1,\mathbf{n},\mathbf{j})^{2} \bigg] \longrightarrow \sigma^2_h,
\end{equation}
\begin{equation}\label{c4}
n^{d}\sum\limits_{\mathbf{j} \in \mathcal{J}} \mathbb{E}\bigg[ U(1,\mathbf{n},\mathbf{j})^{2}\textrm{\large \textbf{1}}_{\{U(1,\mathbf{n},\mathbf{j}) > \varepsilon \sigma_h n^{-d/2} \}}\bigg] \longrightarrow  0\,\,\,\textrm{ for any }  \varepsilon > 0,
\end{equation}
 as $n\rightarrow +\infty$. Indeed, the term $\sum_{\ell=2}^{2^{d}}n^{d/2} \mathscr{S}(\mathbf{n},\ell)$ is asymptotically negligible by \eqref{c2}, the r.v's $n^{d/2}U(1,\mathbf{n},\mathbf{j}) $ are asymptotically independent by \eqref{cd1}, and from \eqref{c3} and the Lindeberg-Feller condition \eqref{c4} we get the convergence in distribution of   $n^{d/2} \mathscr{S}(\mathbf{n},1)$ to $N(0,\sigma^2_h)$  as $n\rightarrow +\infty$.

\bigskip

\noindent Proof of \eqref{cd1}:
Arguing as in \cite{tran90} (see p. 46), we get 
	$
	Q_{1}\leqslant K_1\, r^dp^{d}_{1}\sum_{s=1}^{+\infty}s^{d-1}\varphi(sp_{2})$, where   $K_1>0$.  Since $\varphi(t)\leqslant C\,t^{-\theta}$ and $rp_1\leqslant n$, it follows
\[
Q_{1}\leqslant K_1\,C\,n^{d}p_{2}^{-\theta}\sum_{s=1}^{+\infty}s^{d-1-\theta}\leqslant K_2\,\bigg(\sum_{s=1}^{+\infty}s^{d-1-\theta}\bigg)n^{d-\theta/4}
\]
which yields \eqref{cd1} since $\theta>4d$.

\bigskip

\noindent Proof of \eqref{c2}: We have 
$$n^{d}\mathbb{E}\bigg[\Big( \sum\limits_{\ell=2}^{2^{d}}   \mathscr{S}(\mathbf{n},\ell)\Big)^{2}\bigg]= n^{d}\sum\limits_{\ell=2}^{2^{d}}\mathbb{E}\big[    \mathscr{S}(\mathbf{n},\ell)^{2}\big] + n^{d}\sum_{\ell=2}^{2^d}\sum_{\stackrel{ j=2}{j\neq \ell}}^{2^d}\mathbb{E}\big[  \mathscr{S}(\mathbf{n},\ell)\,\mathscr{S}(\mathbf{n},j)\big].$$
By Cauchy-Schwartz inequality, we get the inequality
$$n^{d}\mathbb{E}\bigg[\Big( \sum\limits_{\ell=2}^{2^{d}}   \mathscr{S}(\mathbf{n},\ell)\Big)^{2}\bigg]\leqslant  n^{d}\sum\limits_{\ell=2}^{2^{d}}\mathbb{E}\big[   \mathscr{S}(\mathbf{n},\ell)^{2}\big] + n^{d}\sum_{\ell=2}^{2^d}\sum_{\stackrel{ j=2}{j\neq \ell}}^{2^d}\mathbb{E}\big[  \mathscr{S}(\mathbf{n},\ell)^2\big]^{1/2}\mathbb{E}\big[  \mathscr{S}(\mathbf{n},j)^2\big]^{1/2}$$
which  shows that it suffices to prove that, for any $\ell\in\{2,\cdots,2^d\}$,  $n^{d}\mathbb{E}\big[   \mathscr{S}(\mathbf{n},\ell)^{2}\big] \longrightarrow 0 $  as $n\rightarrow +\infty$. We will prove it only for $\ell=2$ since the case where $\ell\neq 2$ is similar. Enumerating the $U(2,\mathbf{n},\mathbf{j})$s in an arbitrary way gives $\hat{W}_{1},...,\hat{W}_{M}$, where $M=r^d$. Then $ \mathscr{S}(\mathbf{n},2) = \sum\limits_{\mathbf{j} \in \mathcal{J}}^{}U(2,\mathbf{n},\mathbf{j})=\sum_{j=1}^M\hat{W}_{j}$ and
\begin{eqnarray}\label{s2}
	n^{d}\mathbb{E}\big[ \mathscr{S}(\mathbf{n},2)^{2}\big] & = & n^{d}\sum\limits_{i=1}^{M}var(\hat{W}_{i}) +  n^{d}\sum_{i=1}^M\sum_{\stackrel{ j=1}{j\neq i}}^{M}Cov(\hat{W}_{i},\hat{W}_{j}):= A_{1} + A_{2}.
\end{eqnarray}
As in \cite{dabo16} (see p. 456),  we have
\begin{eqnarray}\label{dvar}
	Var(\hat{W}_{i}) 
	& =&  \sum_{\stackrel{i_{k}=1}{k=1,...,d-1}}^{p_{1}}\,\,\sum_{i_{d}=1}^{p_{2}}Var(\Delta_{\mathbf{i}})
+\sum_{\stackrel{i_{k}=1}{k=1,...,d-1}}^{p_{1}}\,\,\sum_{i_{d}=1}^{p_{2}}\sum_{\stackrel{j_{k}=1}{ \stackrel{ k=1,...,d-1}{ \mathbf{i}\neq \mathbf{j}    }}      }^{p_{1}}\,\,\sum_{j_{d}=1}^{p_{2}}\mathbb{E}\big(\Delta_{\mathbf{i}}\,\Delta_{\mathbf{j}}\big)\nonumber\\
&:=&B_{1}+B_{2}.
\end{eqnarray}
Thanks to the strict stationarity property, we have
\begin{eqnarray}\label{ib1}
B_1&=&p_{1}^{d-1}p_{2}Var(\Delta_{\mathbf{i}})=p_{1}^{d-1}p_{2}n^{-2d}Var\big(h( Z)\big),
\end{eqnarray}
On the other hand, using Lemma A.2 in \cite{dabo16}, we have  
\begin{eqnarray}\label{idelta}
\Big\vert\mathbb{E}\big(\Delta_{\mathbf{i}}\,\Delta_{\mathbf{j}}\big)\Big\vert&=&n^{-2d}\,\Big\vert Cov\left(h(Z_{\mathbf{i}}),h(Z_{\mathbf{j}})\right)\Big\vert\nonumber\\
&\leqslant&  K_2\,n^{-2d}\,
\mathbb{E}\Big(\left\vert h(Z_{\mathbf{i}})\right\vert^4\Big)^{1/4}\mathbb{E}\Big(\left\vert h(Z_{\mathbf{j}})\right\vert^4\Big)^{1/4}\varphi(\parallel\mathbf{i}-\mathbf{j}\parallel_{2})^{1/2}\nonumber\\
&\leqslant&  K_3\,n^{-2d}\ \varphi(\parallel\mathbf{i}-\mathbf{j}\parallel_{2})^{1/2}.
\end{eqnarray}
 Hence
	\begin{eqnarray*}
|B_{2}\vert
	& \leqslant &K_3 n^{-2d}\sum_{\stackrel{i_{k}=1}{k=1,...,d-1}}^{p_{1}}\,\,\sum_{i_{d}=1}^{p_{2}}\sum_{\stackrel{j_{k}=1}{ \stackrel{ k=1,...,d-1}{ \mathbf{i}\neq \mathbf{j}    }}      }^{p_{1}}\,\,\sum_{j_{d}=1}^{p_{2}}\varphi(\parallel\mathbf{i}-\mathbf{j}\parallel_{2})^{1/2}\\
& \leqslant &K_3 n^{-2d}p_{1}^{d-1}p_{2}\sum_{\stackrel{i_{k}=1}{k=1,...,d-1}}^{p_{1}}\,\,\sum_{i_{d}=1}^{p_{2}}\varphi(\parallel\mathbf{i}\parallel_{2})^{1/2}\\
& \leqslant &K_ 4n^{-2d}p_{1}^{d-1}p_2\,\, \sum_{t=1}^{+\infty}t^{d-1}\varphi(t)^{1/2}  \\
& \leqslant &K_4\,C\,n^{-2d}p_{1}^{d-1}p_{2} \sum\limits_{t=1}^{+\infty}t^{d-1-\theta/2} ,
\end{eqnarray*}
that is $\vert B_2\vert\leqslant K_5\,n^{-2d}p_{1}^{d-1}p_{2}$  since $\theta>2d$. Then using this later inequality together with \eqref{ib1} and \eqref{dvar}, we get $Var\left(\hat{W}_i\right)\leqslant  K_6\,n^{-2d}p_{1}^{d-1}p_{2}$. Therefore, since   $M=n^{d}(p_{1} + p_{2})^{-d}\leqslant n^{d}p_{1}^{-d}$, we obtain the inequality
$A_{1}\leqslant  K_6M\,\,n^{-d}p_{1}^{d-1}p_{2}\leqslant  K_6 p_{1}^{-1}p_{2}$ from what we deduce that $A_1\rightarrow 0$ as $n\rightarrow +\infty$ since  $p_{1}^{-1}p_{2}=O(n^{-1/12})$.  For  dealing  with $A_2$ let us first notice that,  as in \cite{dabo16} (see p. 456), we have 
\begin{eqnarray*}
	\vert A_{2} \vert & \leqslant &  n^{d} \underset{ \lVert \mathbf{i}-\mathbf{j} \lVert > p_{2} }{ \sum\limits_{\stackrel{j_{k=1}}{k=1,...,d}}^{n}\sum\limits_{\stackrel{i_{k=1}}{k=1,...,d}}^{n}}\Big\vert \mathbb{E}\big(\Delta_{\mathbf{i}}\Delta_{\mathbf{j}}\big)\Big\vert.
\end{eqnarray*}
Then, using \eqref{idelta} we get the inequality
\begin{eqnarray*}
\vert A_{2} \vert 
	& \leqslant & K_3n^{-d}\sum\limits_{ \lVert\mathbf{i}-\mathbf{j} \lVert > p_{2}} \varphi(\lVert\mathbf{i}-\mathbf{j} \lVert)^{1/2}\\
	& \leqslant &  K_3n^{-d}\sum\limits_{s=1}^{+\infty}\,\,\,\sum\limits_{(s+1)p_{2} \leqslant \lVert \mathbf{i}-\mathbf{j} \lVert = t < (s+2)p_{2}}\varphi(t)^{1/2}\\
	& \leqslant &K_7n^{-d}\sum\limits_{s=1}^{+\infty}s^{d-1}\varphi((s+1)p_{2})^{1/2}\\
	& \leqslant & K_7n^{-d}\sum\limits_{s=1}^{+\infty}s^{d-1}\varphi(s)^{1/2}  \\
	& \leqslant & K_7\,C\,n^{-d}\sum\limits_{s=1}^{+\infty}s^{d-1-\theta/2}
\end{eqnarray*}  
which shows that  $A_2\rightarrow 0$ as $n\rightarrow +\infty$ since $\theta>2d$, and \eqref{c2} follows from  \eqref{s2}.
\bigskip

\noindent Proof of \eqref{c3}:  We have
$ n^{d/2}S_{\mathbf{n}} =  n^{d/2}S_{\mathbf{n}}^{\prime} + n^{d/2}S_{\mathbf{n}}^{\prime\prime},$ where $S_{\mathbf{n}}^{\prime}=\mathscr{S}(\mathbf{n},1)$ and $S_{\mathbf{n}}^{\prime\prime}= \sum_{\ell=2}^{2^{d}} \mathscr{S}(\mathbf{n},\ell) $, hence
\begin{eqnarray}\label{vsprime2}
n^{d}Var\big(S_{\mathbf{n}}^{\prime} \big)&=&n^{d}Var\big(S_{\mathbf{n}}-S_{\mathbf{n}}^{\prime\prime} \big)=n^{d}Var\big(S_{\mathbf{n}}\big)+n^{d}\mathbb{E}\Big(\big(S_{\mathbf{n}}^{\prime\prime} \big)^2\Big) -2n^{d}Cov\big(S_{\mathbf{n}}, S_{\mathbf{n}}^{\prime\prime}\big)\nonumber\\
&=&Var\Big(n^{d/2}S_{\mathbf{n}}\Big)+ n^{d}\mathbb{E}\bigg[\Big( \sum\limits_{\ell=2}^{2^{d}}   \mathscr{S}(\mathbf{n},\ell)\Big)^{2}\bigg] -2 n^{d}\mathbb{E}\Big( S_{\mathbf{n}}\,S_{\mathbf{n}}^{\prime\prime}\Big).
\end{eqnarray}
On the one hand, 
\begin{eqnarray*}
	Var\Big(n^{d/2}S_{\mathbf{n}}\Big) &= & Var\bigg(\dfrac{1}{n^{d/2}}\sum\limits_{\mathbf{i} \in \mathcal{I}_{\mathbf{n}}}\Big(h\left(Z_\mathbf{i}\right)- \mathbb{E}(h(Z_\mathbf{i}))\bigg)\nonumber\\
&=&\dfrac{1}{n^{d}}\sum_{   (\mathbf{i},\mathbf{j} )\in \mathcal{I}_{\mathbf{n}}^2                  }\mathbb{E}\Big[\Big(h\left(Z_\mathbf{i}\right)- \mathbb{E}(h(Z_\mathbf{i})) \Big) \Big( h\left(Z_\mathbf{j}\right)- \mathbb{E}(h(Z_\mathbf{j}))\Big)\Big]\nonumber\\
&=&\dfrac{1}{n^{d}}\sum\limits_{\mathbf{i} \in \mathcal{I}_{\mathbf{n}}}Var(h\left(Z_\mathbf{i}\right))
	+\dfrac{1}{n^{d}}\sum_{    \stackrel{(\mathbf{i},\mathbf{j} )\in \mathcal{I}_{\mathbf{n}}^2}{  \mathbf{i}\neq \mathbf{j}  }                     }Cov(h\left(Z_\mathbf{i}\right),h\left(Z_\mathbf{j}\right)).
\end{eqnarray*}
Thanks to the strict stationnarity, we get 
\[
Var\Big(n^{d/2}S_{\mathbf{n}}\Big)=Var(h(Z))+\frac{1}{n^d}\sum_{    \stackrel{(\mathbf{i},\mathbf{j} )\in \mathcal{I}_{\mathbf{n}}^2}{  \mathbf{i}\neq \mathbf{j}  }                     }Cov(h\left(Z_\mathbf{i}\right),h\left(Z_\mathbf{j}\right)), 
\]
and from  \eqref{condvar} we deduce that $\lim_{n\rightarrow +\infty} Var\Big(n^{d/2}S_{\mathbf{n}}\Big)=\sigma^2_h$. On the other hand, 
\begin{eqnarray*}
\Big\vert n^{d}\mathbb{E}\Big( S_{\mathbf{n}}\,S_{\mathbf{n}}^{\prime\prime}\Big)\Big\vert\leqslant \mathbb{E}\Big( n^{d} S_{\mathbf{n}}^2\Big)^{1/2}\,\mathbb{E}\Big( n^{d}\big(S_{\mathbf{n}}^{\prime\prime}\big)^2\Big)^{1/2}=Var\Big(n^{d/2}S_{\mathbf{n}}\Big)^{1/2} \bigg(  n^{d}\mathbb{E}\bigg[\Big( \sum\limits_{\ell=2}^{2^{d}}   \mathscr{S}(\mathbf{n},\ell)\Big)^{2}\bigg]\bigg)^{1/2},
\end{eqnarray*}
and it follows from  \eqref{c2} that 
$
\lim_{n\rightarrow +\infty  }\bigg(n^{d}\mathbb{E}\Big( S_{\mathbf{n}}\,S_{\mathbf{n}}^{\prime\prime}\Big)\bigg)
=0
$
and, therefore, from \eqref{vsprime2} that
$
\lim_{n\rightarrow +\infty  }\big(n^{d}\mathbb{E}\big(\big(S_{\mathbf{n}}^{\prime} \big)^2\big)\big)
=\lim_{n\rightarrow +\infty} Var\Big(n^{d/2}S_{\mathbf{n}}\Big)=\sigma^2_h$.
Since
\begin{eqnarray*}
	n^{d}\mathbb{E}\Big(\big(S_{\mathbf{n}}^{\prime} \big)^2\Big) & = & n^{d} \sum_{\mathbf{j} \in \mathcal{J}}\mathbb{E}\Big[U(1,\mathbf{n},\mathbf{j})^2\Big]+n^{d} \sum_{\mathbf{j} \in \mathcal{J}} \sum_{\stackrel{\mathbf{l} \in \mathcal{J}}{ \mathbf{l}\neq \mathbf{j}   }}\mathbb{E}\Big[U(1,\mathbf{n},\mathbf{l})\,U(1,\mathbf{n},\mathbf{j})\Big],
\end{eqnarray*}
we get \eqref{c3} if we prove that $n^{d} \sum_{\mathbf{j} \in \mathcal{J}} \sum_{\stackrel{\mathbf{l} \in \mathcal{J}}{ \mathbf{l}\neq \mathbf{j}   }}\mathbb{E}\Big[U(1,\mathbf{n},\mathbf{l})\,U(1,\mathbf{n},\mathbf{j})\Big]\rightarrow 0$ as  $n\rightarrow +\infty$, what is easily obtained since from \eqref{idelta} we have
\begin{eqnarray*}
& &\Big\vert n^{d} \sum_{\mathbf{j} \in \mathcal{J}} \sum_{\stackrel{\mathbf{l} \in \mathcal{J}}{ \mathbf{l}\neq \mathbf{j}   }}\mathbb{E}\Big[U(1,\mathbf{n},\mathbf{l})\,U(1,\mathbf{n},\mathbf{j})\Big] \Big\vert\\
	& \leqslant & n^{d}\sum_{\mathbf{j} \in \mathcal{J}} \sum_{\stackrel{\mathbf{l} \in \mathcal{J}}{ \mathbf{l}\neq \mathbf{j}   }}\sum_{\stackrel{i_{k}=j_{k}(p_{1}+p_{2})+1}{k=1,...,d}}^{j_{k}(p_{1}+p_{2})+p_{1}} \sum_{\stackrel{m_{s}=l_{s}(p_{1}+p_{2})+1}{s=1,...,d}}^{l_{s}(p_{1}+p_{2})+p_{1}}\Big\vert\mathbb{E}\big(\Delta_{\mathbf{i}}\,\Delta_{\mathbf{m}}\big)\Big\vert\\
& \leqslant &K_3 n^{-d}\sum_{\mathbf{j} \in \mathcal{J}} \sum_{\stackrel{\mathbf{l} \in \mathcal{J}}{ \mathbf{l}\neq \mathbf{j}   }}\sum_{\stackrel{i_{k}=j_{k}(p_{1}+p_{2})+1}{k=1,...,d}}^{j_{k}(p_{1}+p_{2})+p_{1}} \sum_{\stackrel{m_{s}=l_{s}(p_{1}+p_{2})+1}{s=1,...,d}}^{l_{s}(p_{1}+p_{2})+p_{1}}\varphi(\parallel\mathbf{i}-\mathbf{m}\parallel_{2})^{1/2}\\
& \leqslant &K_3 n^{-d} \sum_{ \mathbf{i}/\,\Vert \mathbf{i}\Vert_2>0   }\varphi(\parallel\mathbf{i}\parallel_{2})^{1/2}\\
& \leqslant &K_4 n^{-d}\,\, \sum_{t=1}^{+\infty}t^{d-1}\varphi(t)^{1/2}  \\
& \leqslant &K_4\,C\, n^{-d}\sum\limits_{t=1}^{+\infty}t^{d-1-\theta/2} .
\end{eqnarray*}

\bigskip

\noindent Proof of \eqref{c4}: Clearly,    $n^{d/2}\vert \Delta_{\mathbf{i}}\vert \leqslant 2Kn^{-d/2}$, where $K$ is a real that bounds $h$. Hence, since  $p_{1}=\lfloor n^{1/3} \rfloor$,  it follows
\begin{equation}\label{inegu1}
n^{d/2}\vert  U(1,\mathbf{n},\mathbf{j})\vert\leqslant\sum\limits_{\stackrel{i_{k}=j_{k}(p_{1}+p_{2})+1}{k=1,...,d}}^{j_{k}(p_{1}+p_{2})+p_{1}}n^{d/2}\vert \Delta_{\mathbf{i}}\vert \leqslant \dfrac{2K p_1^{d}}{n^{d/2}}\leqslant\dfrac{2K}{n^{d/6}}.
\end{equation}
An  use of Cauchy-Schwartz inequality yields
\begin{eqnarray*}
	& &n^{d}\sum\limits_{\mathbf{j} \in \mathcal{J}} \mathbb{E}\bigg[ U(1,\mathbf{n},\mathbf{j})^{2}\textrm{\large \textbf{1}}_{\{U(1,\mathbf{n},\mathbf{j}) > \varepsilon \sigma_h n^{-d/2} \}}\bigg] \\
&\leqslant& n^{d}\sum\limits_{\mathbf{j} \in \mathcal{T}}^{}\mathbb{E}\big[ U(1,\mathbf{n},\mathbf{j})^{4}]^{1/2}  P\Big( U(1,\mathbf{n},\mathbf{j}) > \varepsilon \sigma_h n^{-d/2}\Big)^{1/2}\\
	&\leqslant& 4K^2n^{-d/3}\sum\limits_{\mathbf{j} \in \mathcal{T}} P\Big( U(1,\mathbf{n},\mathbf{j}) > \varepsilon \sigma_h n^{-d/2}\Big)^{1/2};
\end{eqnarray*}
and since, from \eqref{inegu1}, we have   $\lim_{n\rightarrow +\infty}\dfrac{n^{d/2}U(1,\mathbf{n},\mathbf{j})}{ \sigma_h}=0$,  we deduce that for  $n$ large enough  $ P\Big( U(1,\mathbf{n},\mathbf{j}) > \varepsilon \sigma_h n^{-d/2}\Big)= 0 $ for any  $\mathbf{j} \in  \mathcal{J}$ and, therefore, that \eqref{c4} holds.

\subsection{Technical  lemmas}
\noindent In this section, we give and prove some lemmas that will be useful for proving Theorem \ref{convest}. Note that any linear map   $T \in \mathscr{L}(\mathbb{R}^{p+q},\mathbb{R}^{p+q} )$ can be writen as
\begin{equation*}
	T = \begin{pmatrix}
		T_{11}	&T_{12} \\
		T_{21}&T_{22}
	\end{pmatrix}
\end{equation*} 
with $T_{11} \in \mathscr{L}(\mathbb{R}^{p})$, $T_{22} \in \mathscr{L}(\mathbb{R}^{q})$, $T_{12} \in \mathscr{L}(\mathbb{R}^{q},\mathbb{R}^{p})$  and
$T_{21} \in \mathscr{L}(\mathbb{R}^{p},\mathbb{R}^{q})$. 
Consider then, the following canonical projectors:
\[
	\pi_{1}:\,T\in\mathscr{L}(\mathbb{R}^{p+q}) \mapsto T_{11}  \in\mathscr{L}(\mathbb{R}^{p})\,\,\,\textrm{ and }\,\,\,\pi_{2}:\,T\in\mathscr{L}(\mathbb{R}^{p+q}) \mapsto T_{12}\in  \mathscr{L}(\mathbb{R}^{q},\mathbb{R}^{p}).
\]
So,   the covariance operator $V_{Z} $ of $Z $ can be writen as 

\begin{equation*}
	V_{Z} = \mathbb{E}\bigg(Z\otimes Z\bigg)=
	\begin{pmatrix}
		V_{1}	&V_{12} \\
		V_{21}&V_{2}
	\end{pmatrix},
\end{equation*}
with  $V_{2}=\mathbb{E}\big( Y\otimes Y \big)$ and $V_{21}=V_{12}^{\ast}$, and, putting $\overline{Z}^{(\mathbf{n})}= n^{-d}\displaystyle{\sum_{\mathbf{i} \in \mathcal{I}_{\mathbf{n}}}^{} Z_{\mathbf{i}}}$,  its   empirical estimator  is given by 
	\begin{equation*}
	\widehat{ V}^{(\mathbf{n})}_{Z}  = \dfrac{1}{n^{d}}\displaystyle{\sum_{\mathbf{i} \in \mathcal{I}_{\mathbf{n}}}^{} (Z_{\mathbf{i}}-\overline{Z}^{(\mathbf{n})} )\otimes}(Z_{\mathbf{i}}-\overline{Z}^{(\mathbf{n})})=
	\begin{pmatrix}
		\widehat{V}^{(\mathbf{n})}_{1}	&\widehat{V}^{(\mathbf{n})}_{12} \\
		\widehat{V}^{(\mathbf{n})}_{21}&\widehat{V}^{(\mathbf{n})}_{2}
	\end{pmatrix}   ,
\end{equation*} 
where
\[
\widehat{V}_{2}^{(\mathbf{n})}=\dfrac{1}{n^{d}}\displaystyle{\sum_{\mathbf{i} \in \mathcal{I}_{\mathbf{n}}}^{} (Y _{\mathbf{i}}-\overline{Y }^{(\mathbf{n})})\otimes (Y_{\mathbf{i}}-\overline{Y}^{(\mathbf{n})} )}\,\,\,\textrm{ and }\,\,\,\widehat{V}^{(\mathbf{n})}_{21}=\big(\widehat{V}^{(\mathbf{n})}_{12}\big)^\ast.
\]

\noindent Then, putting
\begin{equation}\label{hn}
\widehat{H}^{(\textbf{n})} = n^{d/2}\Big( 	\widehat{V}^{(\mathbf{n})}_{Z} -V_Z\Big),
\end{equation}
we have $\widehat{H}^{(\textbf{n})}=\widehat{H}^{(\textbf{n})}_1-\widehat{H}^{(\textbf{n})}_2$, where
\begin{equation}\label{h1}
\widehat{H}^{(\textbf{n})}_{1} = \dfrac{1}{n^{d/2}}\sum\limits_{\mathbf{i} \in \mathcal{I}_{\mathbf{n}}}\Big(Z _{\textbf{i}}\otimes Z _{\textbf{i}} - V_Z\Big),
\end{equation}
and  $\widehat{H}^{(\textbf{n})}_{2} = n^{d/2}\overline{Z }^{(\textbf{n})}\otimes\overline{Z }^{(\textbf{n})}$. Let us consider the covariance operator $\Gamma$ of $Z\otimes Z$ given by
\[
\Gamma=\mathbb{E}\Big(\left(Z\otimes Z-V_Z\right)\widetilde{\otimes}\left(Z\otimes Z-V_Z\right)\Big),
\]
where $\widetilde{\otimes}$ denotes the tensor product related to $ \langle\cdot,\cdot\rangle_\mathcal{H}$. Furthermore, considering  an orthonormal basis  $\{\mathcal{U}_k\}_{1\leqslant k\leqslant (p+q)^2}$  of $\mathscr{L}\left(\mathbb{R}^{p+q}\right)$ such that $\mathcal{U}_k$ is an eigenvector of $\Gamma$ 
associated to  the $k$-th largest eigenvalue $\nu_k$, and putting
\[
\Lambda=\sum_{k=1}^{(p+q)^2}\sum_{\ell=1}^{(p+q)^2}\nu_{k}^{1/2}\nu_{\ell}^{1/2}\eta_{k\ell}\,\,\theta_{k\ell}\,\,\mathcal{U}_k\widetilde{\otimes}\mathcal{U}_\ell,
\]
where $\eta_{k\ell}$ and $\theta_{k\ell}$ are defined in \eqref{cvcor2}, we have the following lemma which gives asymptotic normality for the random operator $\widehat{H}^{(\textbf{n})} $.

\bigskip

\begin{lemma}\label{normalite}
 Under the assumptions 2 and  4, if $\varphi(t)=O(t^{-\theta})$ with $\theta>4d$, then   $\widehat{H}^{(\textbf{n})} $ converges in distribution, as $n\rightarrow +\infty$, to a random operator $H$ having a normal distribution in $\mathscr{L}(\mathbb{R}^{p+q})$ with mean equal to $0$ and covariance operator equal to $\Gamma+\Lambda$. 
\end{lemma}

\noindent\textit{Proof}. 
For any $L>0$,   we clearly have $\widehat{H}^{(\mathbf{n})}_{1}=\widehat{H}^{(\mathbf{n})L}_{1}+\widehat{H}^{(\mathbf{n})L\ast}_{1}$, where 
\[
\widehat{H}^{(\textbf{n})L}_{1} = \dfrac{1}{n^{d/2}}\sum\limits_{\mathbf{i} \in \mathcal{I}_{\mathbf{n}}}\Big(Z_\mathbf{i}^{(L)}\otimes Z_\mathbf{i}^{(L)} - V^L_Z\Big)
\,\,\,\textrm{ and }
\,\,\,
\widehat{H}^{(\textbf{n})L\ast}_{1} = \dfrac{1}{n^{d/2}}\sum\limits_{\mathbf{i} \in \mathcal{I}_{\mathbf{n}}}\Big(Z_\mathbf{i}^{(L\ast)}\otimes Z_\mathbf{i}^{(L\ast)} - V^{L\ast}_Z\Big),
\]
with $V^L_Z=\mathbb{E}\left(Z_\mathbf{i}^{(L)}\otimes Z_\mathbf{i}^{(L)}\right)$ and  $V^{L\ast}_Z=\mathbb{E}\left(Z_\mathbf{i}^{(L\ast)}\otimes Z_\mathbf{i}^{(L\ast)}\right)$.   We will  first  get asymptotic normality for $\widehat{H}^{(\textbf{n})L}_{1}$. For any $T\in\mathscr{L}(\mathbb{R}^{p+q})$, we have 
$
\Big<\widehat{H}^{(\textbf{n})L}_{1},T\Big>_{\mathcal{H}}=n^{-d/2}\sum\limits_{\mathbf{i} \in \mathcal{I}_{\mathbf{n}}}\Big(h\left(Z_\mathbf{i}\right)- \mathbb{E}(h(Z_\mathbf{i}))\Big)$,
where  $
h(x)=\Big<x\otimes x,T\Big>_{\mathcal{H}}\textrm{\large \textbf{1}}_{]-\infty,L]}(\Vert x\Vert_{\mathbb{R}^{p+q}})$.
First,
\[
Var\left(h(Z)\right)=Var\left( \Big<Z_\mathbf{i}^{(L)}\otimes Z_\mathbf{i}^{(L)} - V^L_Z,T \Big>_{\mathcal{H}}  \right)=\Big<\Gamma_L(T),T \Big>_{\mathcal{H}} ,
\]
where
\[
\Gamma_L=\mathbb{E}\Big(\left(Z^{(L)}\otimes Z^{(L)}-V^L_Z\right)\widetilde{\otimes}\left(Z^{(L)}\otimes Z^{(L)}-V^L_Z\right)\Big).
\]
On the other hand, considering  an orthonormal basis  $\{\mathcal{U}^L_k\}_{1\leqslant k\leqslant (p+q)^2}$  of $\mathscr{L}\left(\mathbb{R}^{p+q}\right)$ such that $\mathcal{U}^L_k$ is an eigenvector of $\Gamma_L$ 
associated to  the $k$-th largest eigenvalue $\nu^{(L)}_k$,  we have $T=\sum_{k=1}^{(p+q)^2}t_k^L\,\mathcal{U}^L_k$, where   $t^L_{k} =\big<T,\mathcal{U}^L_{k}\big>_{\mathcal{H}}$. Using   Assumption \ref{ass3}-$(ii)$ (see \eqref{cvcor2}), we get
\begin{eqnarray*}
& &\lim_{n\rightarrow +\infty}\bigg\{\dfrac{1}{n^{d}}\sum_{    \stackrel{(\mathbf{i},\mathbf{j} )\in \mathcal{I}_{\mathbf{n}}^2}{  \mathbf{i}\neq \mathbf{j}  } }Cov\Big(h\left(Z_\mathbf{i}\right),h\left(Z_\mathbf{j}\right)\Big)\bigg\}\\
&=&\lim_{n\rightarrow +\infty}\bigg\{ \dfrac{1}{n^{d}}\sum_{    \stackrel{(\mathbf{i},\mathbf{j} )\in \mathcal{I}_{\mathbf{n}}^2}{  \mathbf{i}\neq \mathbf{j}  }                     }\mathbb{E}\Big[\Big<Z _{\textbf{i}}^{(L)}\otimes Z _{\textbf{i}}^{(L)} -V_Z^L,T\Big>_{\mathcal{H}}\Big<Z _{\textbf{j}}^{(L)}\otimes Z _{\textbf{j}}^{(L)} - V_Z^L,T\Big>_{\mathcal{H}}\Big]\bigg\}\\
&=&\lim_{n\rightarrow +\infty}\bigg\{\dfrac{1}{n^{d}}\sum_{k=1}^{(p+q)^2}\sum_{\ell=1}^{(p+q)^2}t^L_{k}t^L_{\ell}\sum_{    \stackrel{(\mathbf{i},\mathbf{j} )\in \mathcal{I}_{\mathbf{n}}^2}{  \mathbf{i}\neq \mathbf{j}  } }\mathbb{E}\Big[\Big<Z^{(L)} _{\textbf{i}}\otimes Z^{(L)} _{\textbf{i}} -V^L_Z,\mathcal{U}^L_{k}\Big>_{\mathcal{H}}\Big<Z ^{(L)}_{\textbf{j}}\otimes Z ^{(L)}_{\textbf{j}} - V^L_Z,\mathcal{U}^L_{\ell}\Big>_{\mathcal{H}}\Big]\bigg\}\\
&=&\sum_{k=1}^{(p+q)^2}\sum_{\ell=1}^{(p+q)^2}t^L_{k}t^L_{\ell}\eta^{(L)}_{k\ell} \left(\nu_k^{(L)}\nu_\ell^{(L)}\right)^{1/2}\lim_{n\rightarrow +\infty}\bigg\{\dfrac{1}{n^{d}}\sum_{    \stackrel{(\mathbf{i},\mathbf{j} )\in \mathcal{I}_{\mathbf{n}}^2}{  \mathbf{i}\neq \mathbf{j}  } }\Theta_{k\ell}^{(L)}(\Vert\mathbf{i}-\mathbf{j}\Vert_2)\bigg\}\\
&=&\sum_{k=1}^{(p+q)^2}\sum_{\ell=1}^{(p+q)^2}t^L_{k}t^L_{\ell}\eta^{(L)}_{k\ell} \left(\nu_k^{(L)}\nu_\ell^{(L)}\right)^{1/2}\theta_{k\ell}^{(L)}.\\
\end{eqnarray*}
Since
\[
t^L_{k}t^L_{\ell}=\big<T,\mathcal{U}^L_{k}\big>_{\mathcal{H}}\big<T,\mathcal{U}^L_{\ell}\big>_{\mathcal{H}}=\Big\langle (\mathcal{U}^L_{k}\widetilde{\otimes}\mathcal{U}^L_{\ell})(T),T\Big\rangle,
\]
we finally obtain
\begin{eqnarray*}
\lim_{n\rightarrow +\infty}\bigg\{\dfrac{1}{n^{d}}\sum_{    \stackrel{(\mathbf{i},\mathbf{j} )\in \mathcal{I}_{\mathbf{n}}^2}{  \mathbf{i}\neq \mathbf{j}  } }Cov\Big(h\left(Z_\mathbf{i}\right),h\left(Z_\mathbf{j}\right)\Big)\bigg\}
&=&\big< \Lambda_L(T),T\big>_{\mathcal{H}},
\end{eqnarray*}
where
\[
\Lambda_L=\sum_{k=1}^{(p+q)^2}\sum_{\ell=1}^{(p+q)^2}\eta^{(L)}_{k\ell} \left(\nu_k^{(L)}\nu_\ell^{(L)}\right)^{1/2}\theta_{k\ell}^{(L)} \mathcal{U}^L_{k}\widetilde{\otimes}\mathcal{U}^L_{\ell}.
\]
Then,  from   Theorem \ref{clt},we deduce  that $\Big<\widehat{H}^{(\mathbf{n})L}_{1},T\Big>_{\mathcal{H}}$ converges in distribution to $\mathcal{N}(0,\sigma^2_L)$, where $\sigma^2_L=\Big<\Gamma_L(T),T \Big>_{\mathcal{H}} +\Big<\Lambda_L(T),T \Big>_{\mathcal{H}} $.
 By  Levy's Theorem, this later property implies
\begin{equation*}\label{c1}
 \mathbb{E}\bigg[\exp\bigg(i\Big<\widehat{H}^{(\mathbf{n})L}_{1},T\Big>_{\mathcal{H}} \bigg)\bigg] \longrightarrow \exp\bigg(-\frac{1}{2}\big<(\Gamma_L+\Lambda_L)(T),T\big>_{\mathcal{H}} \bigg)\,\,\,\textrm{ as }n\rightarrow +\infty,
\end{equation*}
for all $T\in\mathscr{L}(\mathbb{R}^{p+q})$, what is equivalent to  the convergence in distribution of $\widehat{H}^{(\textbf{n})L}_1 $, as $n\rightarrow +\infty$, to a random operator $H_L$ having a normal distribution in $\mathscr{L}(\mathbb{R}^{p+q})$ with mean equal to $0$ and covariance operator equal to $\Gamma_L+\Lambda_L$. Secondly, we will prove 
that
\begin{equation*} 
 \mathbb{E}\bigg[\exp\bigg(i\Big<\widehat{H}^{(\mathbf{n})}_{1},T\Big>_{\mathcal{H}} \bigg)\bigg] \longrightarrow \exp\bigg(-\frac{1}{2}\big<(\Gamma+\Lambda)(T),T\big>_{\mathcal{H}} \bigg) \,\,\,\textrm{ as }n\rightarrow +\infty,
\end{equation*}
for all $T\in\mathscr{L}(\mathbb{R}^{p+q})$. Putting $\sigma^2=\big<(\Gamma+\Lambda)(T),T\big>_{\mathcal{H}}$, we have
\begin{eqnarray*}
	& & \Big|  \mathbb{E}\bigg[\exp\bigg(i\Big<\widehat{H}^{(\mathbf{n})}_{1},T\Big>_{\mathcal{H}} \bigg)\bigg]  - \exp\Big(-\dfrac{\sigma^2}{2} \Big) \Big|\\
	&=&\Big\vert \mathbb{E}\bigg[\exp\bigg(i\Big<\widehat{H}^{(\mathbf{n})L}_{1},T\Big>_{\mathcal{H}} \bigg)\times \exp\bigg(i\Big<\widehat{H}^{(\mathbf{n})L\ast}_{1},T\Big>_{\mathcal{H}} \bigg)\bigg]  - \exp\Big(-\dfrac{\sigma^2}{2} \Big) \Big| \\
	& \leqslant&\Big| \mathbb{E}\bigg[\exp\bigg(i\Big<\widehat{H}^{(\mathbf{n})L}_{1},T\Big>_{\mathcal{H}}\Big) \bigg] - \exp\Big( -\dfrac{ \sigma^2_L}{2}\Big) \Big| + \Big|\exp\Big(-\dfrac{\sigma^2_L}{2} \Big)- \exp\Big(-\dfrac{\sigma^2}{2} \Big) \Big|\\
	& +& \Big|\mathbb{E}\bigg[\exp\bigg(i\Big<\widehat{H}^{(\mathbf{n})L\ast}_{1},T\Big>_{\mathcal{H}}\Big) \bigg]-1 \Big|.
\end{eqnarray*}
Since $\widehat{H}^{(\mathbf{n})L}_{1}$ converges in distribution to $N(0,\Gamma_L+\Lambda_L)$ as $n\rightarrow +\infty$, we have $\Big| \mathbb{E}\bigg[\exp\bigg(i\Big<\widehat{H}^{(\mathbf{n})L}_{1},T\Big>_{\mathcal{H}}\Big) \bigg] - \exp\Big( -\dfrac{ \sigma^2_L}{2}\Big) \Big|\rightarrow 0$. Further, from dominated convergence theorem, we obtain that $\Gamma_L\rightarrow \Gamma$, as $L\rightarrow +\infty$. Then, using Lemma 1 in \cite{ferre2003} we can deduce that $\nu_k^{(L)}\rightarrow \nu_k$ and $\mathcal{U}^L_k\rightarrow\mathcal{U}_k$ as $L\rightarrow +\infty$. These  later properties together with \eqref{cvcor2} imply that $\Lambda_L\rightarrow \Lambda$, as $L\rightarrow +\infty$. Hence   $\sigma_L^2\rightarrow \sigma^2$   and  $\Big|\exp\Big(-\dfrac{\sigma^2_L}{2} \Big)- \exp\Big(-\dfrac{\sigma^2}{2} \Big) \Big|\rightarrow 0$ as $L\rightarrow +\infty$. It remains to prove that $ \Big|\mathbb{E}\bigg[\exp\bigg(i\Big<\widehat{H}^{(\mathbf{n})L\ast}_{1},T\Big>_{\mathcal{H}}\Big) \bigg]-1 \Big|\rightarrow 0$ as $L\rightarrow +\infty$.   Clearly, $\mathbb{E}\Big[\Big<\widehat{H}^{(\mathbf{n})L\ast}_{1},T\Big>_{\mathcal{H}}^2 \Big] =A+B$, where $
	A=  n^{-d}\sum\limits_{\mathbf{i} \in \mathcal{I}_{\mathbf{n}}} \mathbb{E}\Big[ \Big<Z_{\mathbf{i}}^{(L\ast)}\otimes Z_{\mathbf{i}}^{(L\ast)} -V_{Z}^{L\ast}, T\big>_{\mathcal{H}}^{2}\Big]$ and
\begin{eqnarray*}
	B=\dfrac{1}{n^{d}}\sum_{    \stackrel{(\mathbf{i},\mathbf{j} )\in \mathcal{I}_{\mathbf{n}}^2}{  \mathbf{i}\neq \mathbf{j}  }                     } \mathbb{E}\Big[\Big<Z_{\mathbf{i}}^{(L\ast)}\otimes Z_{\mathbf{i}}^{(L\ast)} -V_{Z}^{L\ast}, T\big>_{\mathcal{H}}\Big<Z_{\mathbf{j}}^{(L\ast)}\otimes Z_{\mathbf{j}}^{(L\ast)} -V_{Z}^{L\ast}, T\big>_{\mathcal{H}}\Big].
\end{eqnarray*} 
By Cauchy-Schwartz inequality and the strict stationarity,
\begin{eqnarray}\label{inegA}
A&\leqslant& \Big\Vert  T\Big\Vert _{\mathcal{H}}^{2} \mathbb{E}\Big[ \Big\Vert Z_{\mathbf{i}}^{(L\ast)}\otimes Z_{\mathbf{i}}^{(L\ast)} -V_{Z}^{L\ast}\Big\Vert _{\mathcal{H}}^2\Big]
\leqslant 2\Big\Vert  T\Big\Vert _{\mathcal{H}}^{2} \mathbb{E}\Big[ \Big\Vert Z_{\mathbf{i}}^{(L\ast)}\otimes Z_{\mathbf{i}}^{(L\ast)} \Big\Vert _{\mathcal{H}}^2+\Big\Vert  V^{L\ast}_Z\Big\Vert _{\mathcal{H}}^2\Big]\nonumber\\
&\leqslant &4\Big\Vert  T\Big\Vert _{\mathcal{H}}^{2} \mathbb{E}\Big[ \left\Vert Z_{\mathbf{i}}^{(L\ast)}\right\Vert^4_{\mathbb{R}^{p+q}}\Big].
\end{eqnarray}
Using again Cauchy-Schwartz inequality followed by Markov inequality we obtain
\begin{eqnarray}\label{moment4}
 \mathbb{E}\Big[ \left\Vert Z_{\mathbf{i}}^{(L\ast)}\right\Vert^4_{\mathbb{R}^{p+q}}\Big]&= &\mathbb{E}\Big[ \left\Vert Z_\mathbf{i}\right\Vert^4_{\mathbb{R}^{p+q}}\textrm{\large \textbf{1}}_{\{ \left\Vert  Z_\mathbf{i}\right\Vert_{\mathbb{R}^{p+q}} > L \}} \Big]
\leqslant \mathbb{E}\Big[ \left\Vert Z_\mathbf{i}\right\Vert^8_{\mathbb{R}^{p+q}} \Big] ^{1/2}P\Big(\left\Vert  Z_\mathbf{i}\right\Vert_{\mathbb{R}^{p+q}} > L\Big) ^{1/2}\nonumber\\
&\leqslant &\frac{ \mathbb{E}\Big[ \left\Vert Z_\mathbf{i}\right\Vert^8_{\mathbb{R}^{p+q}} \Big] ^{1/2}\mathbb{E}\Big[ \left\Vert Z_\mathbf{i}\right\Vert_{\mathbb{R}^{p+q}} \Big] ^{1/2} }{ L^{1/2} },
\end{eqnarray}
and it follows from \eqref{inegA} and  \eqref{moment4} that $A\rightarrow 0$ as $L\rightarrow +\infty$. On the other hand, considering  an orthonormal basis  $\{\mathcal{U}^{L\ast}_k\}_{1\leqslant k\leqslant (p+q)^2}$  of $\mathscr{L}\left(\mathbb{R}^{p+q}\right)$ such that $\mathcal{U}^{L\ast}_k$ is an eigenvector of 
\[
\Gamma_{L\ast}=\mathbb{E}\Big(\left(Z^{(L\ast)}\otimes Z^{(L\ast)}-V^{L\ast}_Z\right)\widetilde{\otimes}\left(Z^{(L\ast)}\otimes Z^{(L\ast)}-V^{L\ast}_Z\right)\Big).
\]
associated to  the $k$-th largest eigenvalue $\nu^{(L\ast)}_k$,  we have $T=\sum_{k=1}^{(p+q)^2}t_k^{L\ast}\,\mathcal{U}^{L\ast}_k$, where   $t_k^{L\ast} =\big<T,\mathcal{U}^{L\ast}_k\big>_{\mathcal{H}}$. Using   Assumption \ref{ass3}-$(ii)$, we get
\begin{eqnarray}\label{inegb}
	B & = &\dfrac{1}{n^{d}}\sum\limits_{\stackrel{(\mathbf{i,j}) \in \mathcal{I}_{\mathbf{n}}^{2}}{\mathbf{i}\neq \mathbf{j}    }} \mathbb{E}\Big[\big<Z_{\mathbf{i}}^{(L\ast)}\otimes Z_{\mathbf{i}}^{(L\ast)} -V_{Z}^{L\ast}, T\big>_{\mathcal{H}}\big<Z_{\mathbf{j}}^{(L\ast)}\otimes Z_{\mathbf{j}}^{(L\ast)} -V_{Z}^{L\ast}, T\big>_{\mathcal{H}}\Big]\nonumber\\
& = &\dfrac{1}{n^{d}}\sum_{k=1}^{(p+q)^2}\sum_{\ell=1}^{(p+q)^2}t_k^{L\ast}t_\ell^{L\ast}\sum\limits_{\stackrel{(\mathbf{i,j}) \in \mathcal{I}_{\mathbf{n}}^{2}}{\mathbf{i}\neq \mathbf{j}    }} \mathbb{E}\Big[\big<Z_{\mathbf{i}}^{(L\ast)}\otimes Z_{\mathbf{i}}^{(L\ast)} -V_{Z}^{L\ast},\mathcal{U}^{L\ast}_k\big>_{\mathcal{H}}\big<Z_{\mathbf{j}}^{(L\ast)}\otimes Z_{\mathbf{j}}^{(L\ast)} -V_{Z}^{L\ast}, \mathcal{U}^{L\ast}_\ell\big>_{\mathcal{H}}\Big]\nonumber\\
& = &\dfrac{1}{n^{d}}\sum_{k=1}^{(p+q)^2}\sum_{\ell=1}^{(p+q)^2}t_k^{L\ast}t_\ell^{L\ast}\eta^{(L\ast)}_{k\ell} \left(\nu_k^{(L\ast)}\nu_\ell^{(L\ast)}\right)^{1/2}\sum\limits_{ 0<\parallel \mathbf{i}-\mathbf{j} \parallel_{2}}\Theta^{(L\ast)}_{k\ell}(\lVert\mathbf{i}-\mathbf{j} \lVert_{2}),
\end{eqnarray}
where   $ \nu^{(L\ast)}_k= \mathbb{E}\Big[ \big<Z_{\mathbf{i}}^{(L\ast)}\otimes Z_{\mathbf{i}}^{(L\ast)} -V_{Z}^{L\ast},\mathcal{U}^{L\ast}_k\big>_{\mathcal{H}}^{2}\Big]$. Since 
\[
\dfrac{1}{n^{d}}\sum\limits_{ 0<\parallel \mathbf{i}-\mathbf{j} \parallel_{2}}\Theta^{(L\ast)}_{k\ell}(\lVert\mathbf{i}-\mathbf{j} \lVert_{2})\leqslant\sum_{t=1}^{+\infty}t^{d-1}\Theta^{(L\ast)}_{k\ell}(t)<+\infty,
\]
it follows from the inequalities 
\[
\left\vert  t_k^{L\ast}t_\ell^{L\ast}\right\vert\leqslant \Vert T\Vert_\mathcal{H}^2\Vert \mathcal{U}^{L\ast}_k\Vert_\mathcal{H}\Vert \mathcal{U}^{L\ast}_\ell\Vert_\mathcal{H}=\Vert T\Vert_\mathcal{H}^2,\,\,\, \Big(\nu^{(L\ast)}_k\nu^{(L\ast)}_\ell\Big)^{1/2}\leqslant \nu^{(L\ast)}_k+\nu^{(L\ast)}_\ell,
\]
\begin{eqnarray*}
	\nu^{(L\ast)}_k&\leqslant &  \mathbb{E}\Big[ \big\Vert Z_{\mathbf{i}}^{(L\ast)}\otimes Z_{\mathbf{i}}^{(L\ast)} -V_{Z}^{L\ast}\Vert_{\mathcal{H}}^{2}\Vert  \mathcal{U}^{L\ast}_k\big\Vert_{\mathcal{H}}^{2}\Big]\\
&\leqslant&   2\mathbb{E}\Big[ \big\Vert Z_{\mathbf{i}}^{(L\ast)}\otimes Z_{\mathbf{i}}^{(L\ast)}\Vert_{\mathcal{H}}^{2} +\big\Vert V_{Z}^{L*}\Vert_{\mathcal{H}}^{2}\Big]\\
&\leqslant&  4\mathbb{E}\Big[ \big\Vert Z_{\mathbf{i}}^{(L\ast)}\big\Vert^{4} _{\mathbb{R}^{p+q}} \Big],
\end{eqnarray*}
and the equations \eqref{moment4} and \eqref{inegb} that $\vert B\vert\leqslant C/L^{1/2}$, where $C$ is a positive constant, what allows to deduce that  $B\rightarrow 0$ as $L\rightarrow +\infty$. Therefore, $\mathbb{E}\Big[\Big<\widehat{H}^{(\mathbf{n})L\ast}_{1},T\Big>_{\mathcal{H}}^2 \Big] \rightarrow 0$ as $L\rightarrow +\infty$, and from Markov inequality, we deduce that $\Big<\widehat{H}^{(\mathbf{n})L\ast}_{1},T\Big>_{\mathcal{H}}$ converges in probability to $0$ as $L\rightarrow +\infty$, what implies that $ \Big|\mathbb{E}\bigg[\exp\bigg(i\Big<\widehat{H}^{(\mathbf{n})L\ast}_{1},T\Big>_{\mathcal{H}}\Big) \bigg]-1 \Big|\rightarrow 0$ as $L\rightarrow +\infty$. Finally, we have proved that $\widehat{H}^{(\textbf{n})}_1$ converges in distribution, as $n\rightarrow +\infty$ to a random variable $H$ having the normal distribution in $\mathscr{L}(\mathbb{R}^{p+q})$ with mean $0$ and covariance operator $\Gamma +\Lambda$. Furthermore, the random operators defined in \eqref{hn} and \eqref{h1} verify  $\widehat{H}^{(\textbf{n})}=\widehat{H}^{(\textbf{n})}_1-\widehat{H}^{(\textbf{n})}_2$, where $\widehat{H}^{(\textbf{n})}_{2} = n^{d/2}\overline{Z }^{(\textbf{n})}\otimes\overline{Z }^{(\textbf{n})}$. Using $
h(x)=\Big<x,z\Big>_{\mathbb{R}^{p+q}}\textrm{\large \textbf{1}}_{]-\infty,L]}(\Vert x\Vert_{\mathbb{R}^{p+q}})$ for any $z\in\mathbb{R}^{p+q}$ and \eqref{cvcor} we obtain,  by a similar reasoning than above, that $n^{d/2}\overline{Z }^{(\textbf{n})}$ converges in distribution, as  $n\rightarrow +\infty$, to a random variable with normal distribution in $\mathbb{R}^{p+q}$. Thus    $\widehat{H}^{(\textbf{n})}_{2} $ converges in probability to $0$, as   $n\rightarrow +\infty$, since  $\widehat{H}^{(\textbf{n})}_{2} =n^{-d/2}\left( n^{d/2}\overline{Z }^{(\textbf{n})}\right)\otimes\left( n^{d/2}\overline{Z }^{(\textbf{n})}\right)$. From Slutsky's theorem we deduce that $\widehat{H}^{(\textbf{n})}$ converges in distribution, as $n\rightarrow +\infty$, to  $H$.  
\hfill$\Box$

\bigskip

\noindent There exist
$r\in I$ and $\left(  m_{1},\cdots,m_{r}\right)  \in I^{r}$ such that
\ $m_{1}+\cdots+m_{r}=p$, and
$\xi_{K_{\tau (  1)  }}    =\cdots=\xi_{K_{\tau(
m_{1})  }}>\xi_{K_{\tau (  m_{1}+1)  }}=\cdots
=\xi_{K_{\tau (  m_{1}+m_{2})  }}>\cdots
\cdots   >\xi_{K_{\tau (  m_{1}+\cdots+m_{r-1}+1)  }}=\cdots
=\xi_{K_{\tau (  m_{1}+\cdots+m_{r})  }}$.
Then considering \ the set $E=\left\{  \ell\in\mathbb{N}^{\ast};1\leq \ell\leq
r,m_{\ell}\geq2\right\}  $ and putting $m_{0}:=0$ and $$F_{\ell}:=\left\{\sum_{k=0}^{\ell-1}m_k+1,\cdots,\sum_{k=0}^{\ell}m_k-1\right\}  $$
\ ($\ell\in\left\{  1,\cdots,r\right\}  $), we  have

\begin{lemma}\label{diffcrit} 	
 Under the assumptions 1 to 4,  suppose that  $\varphi(t)=O(t^{-\theta})$ with $\theta>4d$. If $I \neq \emptyset$ and $\gamma\in]0,1/2[$ , then for all $\ell \in I$ and $i \in F_{\ell}$,  $n^{d\gamma}(\widehat{\xi}_{K_{\tau (i)}}^{(\mathbf{n})} - \widehat{\xi}_{K_{\tau (i+1)}}^{(\mathbf{n})})$ converges in probability to 0  as  $n \rightarrow +\infty.$
\end{lemma}

\noindent\textit{Proof}. First,
\begin{eqnarray}\label{dvlpmt}
&& n^{d/2}\widehat{\xi}_{K}^{\left(  n\right)  }\nonumber\\
  & =&\Big\Vert  n^{d/2}\left(
\widehat{V}_{12}^{\mathbf{(n)}}-V_{12}\right)  - n^{d/2}\left(
\widehat{V}_{1}^{\mathbf{(n)}}-V_{1}\right)  \widehat{\Pi}_{K}^{\left(
n\right)  }\widehat{V}_{12}^{\mathbf{(n)}} \nonumber\\
& & +V_{1}\widehat{\Pi}_{K}^{\left(  n\right)  }\left[   n^{d/2} \left(
\widehat{V}_{1}^{\mathbf{(n)}}left.-V_{1}\right)  \right]  \Pi_{K}\widehat{V}_{12}^{\mathbf{(n)}} -V_{1}\Pi_{K}\left[  n^{d/2}\left(  \widehat{V}_{12}^{\mathbf{(n)}}-V_{12}\right)  \right]  + n^{d/2}\mathcal{D}_{K}\Big\Vert_\mathcal{H} \nonumber \\
  &=&\Big\lVert \widehat{\Phi}^{(\mathbf{n})}_{K}(\widehat{H}^{(\textbf{n})})  + n^{d/2}\mathcal{D}_{K} \Big\lVert_\mathcal{H},
\end{eqnarray}
where $\widehat{\Phi}^{(\mathbf{n})}_{K}$ is the random operator from $\mathscr{L}(\mathbb{R}^{p+q})$ to $\mathscr{L}(\mathbb{R}^{q},\mathbb{R}^{p})$ given by
\[
\widehat{\Phi}^{(\mathbf{n})}_{K}(T) = \pi_{2}(T)-\pi_{1}(T)\widehat{\Pi}_{K}^{\mathbf{(n)}}\widehat{V}_{12}^{\mathbf{(n)}}+V_{1}\widehat{\Pi}_{K}^{\mathbf{(n)}}\pi_{1}(T)\Pi_{K}\widehat{V}_{12}^{\mathbf{(n)}}-V_{1}\Pi_{K}\pi_{2}(T),
\]
and $\mathcal{D}_K=V_{12} - V_{1}\Pi_{K}V_{12}$.  Since, obviously,    the integers  $i$ and  $i+1$ both belong to the set $\left\{  m_{0}+m_{1}+\cdots+m_{\ell-1}+1,\cdots,m_{0}+m_{1}+\cdots
+m_{\ell}\right\}  $, we have $\Vert\mathcal{D}_{K_{\tau (i)}}\Vert_\mathcal{H}=\Vert \mathcal{D}_{K_{\tau (i+1)}}\Vert_\mathcal{H}=\mu_\ell$, where $\mu_\ell=\xi_{K_{\tau (m_{0}+m_{1}+\cdots
+m_{\ell})}}$. If $\mu_\ell=0$, then $\mathcal{D}_{K_{\tau (i)}}= \mathcal{D}_{K_{\tau (i+1)}}=0$ and, therefore,
\begin{eqnarray*}
\left|  n^{d\gamma}\left(  \widehat{\xi}_{K_{\tau\left(  i\right)  }}^{\left(
\mathbf{n}\right)  }-\widehat{\xi}_{K_{\tau\left(  i+1\right)  }}^{\left(  \mathbf{n}\right)
}\right)  \right|  &=&n^{d(\gamma-1/2)}\left(  \left\|  \widehat{\Phi}_{K_{\tau\left(  i\right)  }%
}^{\left(  \mathbf{n}\right)  }\left(  \widehat{H}^{\left(  \mathbf{n}\right)  }\right)
\right\|_\mathcal{H}  -\left\|  \widehat{\Phi}_{K_{\tau\left(  i+1\right)  }}^{\left(
\mathbf{n}\right)  }\left(  \widehat{H}^{\left(  \mathbf{n}\right)  }\right)  \right\| _\mathcal{H}  \right)\\
 & \leq &n^{d(\gamma-1/2)}\left(  \left\|  \left(  \widehat
{\Phi}_{K_{\tau\left(  i\right)  }}^{\left(  \mathbf{n}\right)  }-\widehat{\Phi
}_{K_{\tau\left(  i+1\right)  }}^{\left(  \mathbf{n}\right)  }\right)  \left(
\widehat{H}^{\left(  \mathbf{n}\right)  }\right)  \right\|_\mathcal{H}   \right) \\
& \leq &n^{d(\gamma-1/2)}\bigg(\left\Vert  \widehat{\Phi}_{K_{\tau\left(  i\right)  }%
}^{\left(  \mathbf{n}\right)  }\right\Vert_\infty+\left\Vert\widehat{\Phi}_{K_{\tau\left(  i+1\right)  }%
}^{\left(  \mathbf{n}\right)  }\right\Vert  _\infty\bigg)\left\|  \widehat{H}^{\left(  \mathbf{n}\right)
}\right\| _\mathcal{H}  .
\end{eqnarray*}
Since  $\gamma<1/2$, we deduce from this inequality,  Lemma \ref{normalite} and the almost sure convergences of $\widehat{\Phi}_{K_{\tau\left(  i\right)  }%
}^{\left(  \mathbf{n}\right)  }$ and $\widehat{\Phi}_{K_{\tau\left(  i+1\right)  }%
}^{\left(  \mathbf{n}\right)  }$ to   $\Phi_{K_{\tau\left(  i\right)  }} $ and $\Phi_{K_{\tau\left(  i\right)  }}$ respectively, where
\[
\Phi_{K}(T) = \pi_{2}(T)-\pi_{1}(T)\Pi_{K}V_{12}+V_{1}\Pi_{K}\pi_{1}(T)\Pi_{K}V_{12}-V_{1}\Pi_{K}\pi_{2}(T),
\]
 that $ n^{d\gamma}\left(  \widehat{\xi}_{K_{\tau\left(  i\right)  }}^{\left(
\mathbf{n}\right)  }-\widehat{\xi}_{K_{\gamma\left(  i+1\right)  }}^{\left(  \mathbf{n}\right)
}\right)$ converges in probability to $0$ as $n\rightarrow +\infty$.  If $\mu_\ell\neq 0$, then, arguing as in \cite{mbina2023} (see the proof of Lemma 3), we obtain
\begin{align*}
& \left\vert n^{d\gamma}\left(  \widehat{\xi}_{K_{\tau\left(  i\right)  }}^{\left(
\mathbf{n}\right)  }-\widehat{\xi}_{K_{\tau\left(  i+1\right)  }}^{\left( \mathbf{n}\right)
}\right) \right\vert \\
& \leqslant\frac{ n^{d(\gamma-1)}\bigg(\left\Vert  \widehat{\Phi}_{K_{\tau\left(  i\right)  }%
}^{\left(  \mathbf{n}\right)  }\right\Vert_\infty^2+\left\Vert\widehat{\Phi}_{K_{\tau\left(  i+1\right)  }%
}^{\left(  \mathbf{n}\right)  }\right\Vert  _\infty^2\bigg)\left\|  \widehat{H}^{\left(  \mathbf{n}\right)
}\right\|  _\mathcal{H} }{ \left\|  n^{-d/2}\widehat{\Phi}_{K_{\tau\left(
i\right)  }}^{\left(  \mathbf{n}\right)  }\left(  \widehat{H}^{\left(   \mathbf{n}\right)
}\right)  +\mathcal{D}_{K_{\sigma\left(  i\right)  }}\right\|  _\mathcal{H} +\left\|
n^{-d/2}\widehat{\Phi}_{K_{\tau\left(  i+1\right)  }}^{\left(  \mathbf{n}\right)
}\left(  \widehat{H}^{\left(   \mathbf{n}\right)  }\right)  +\mathcal{D}_{K_{\tau\left(
i+1\right)  }}\right\|  _\mathcal{H} }\\
& +2\mu_\ell \,\frac{n^{d(\gamma-1/2})\bigg(\left\Vert  \widehat{\Phi}_{K_{\tau\left(  i\right)  }%
}^{\left(  \mathbf{n}\right)  }\right\Vert_\infty+\left\Vert\widehat{\Phi}_{K_{\tau\left(  i+1\right)  }%
}^{\left(  \mathbf{n}\right)  }\right\Vert  _\infty\bigg)\left\|  \widehat{H}^{\left(  \mathbf{n}\right)
}\right\|  _\mathcal{H} }{\left\|  n^{-d/2}\widehat{\Phi}_{K_{\tau\left(
i\right)  }}^{\left(  \mathbf{n}\right)  }\left(  \widehat{H}^{\left(   \mathbf{n}\right)
}\right)  +\mathcal{D}_{K_{\tau\left(  i\right)  }}\right\| _\mathcal{H}  +\left\|
n^{-d/2}\widehat{\Phi}_{K_{\tau\left(  i+1\right)  }}^{\left(  \mathbf{n}\right)
}\left(  \widehat{H}^{\left(   \mathbf{n}\right)  }\right)  +\mathcal{D}_{K_{\tau\left(
i+1\right)  }}\right\|  _\mathcal{H} }.
\end{align*}
As above, the numerators converge in probability  to 0  as  $n \rightarrow +\infty$. Then the required result is deduced from the fact that  the denominator converges in probability to $2\mu_\ell$ as  $n \rightarrow +\infty$.

\subsection{Proof of Theorem \ref{convest}}
\noindent Since $\left\{
\widehat{\tau}_ \mathbf{n} =\tau\right\}  \cap\left\{  \widehat
{s}_  \mathbf{n}  =s\right\}  \subset\left\{  \widehat{I}_{1}^{\left(
\mathbf{n}\right)  }=I_{1}\right\}  $, it suffices to prove that 
\begin{equation}\label{conv1}
\lim_{n\rightarrow +\infty}P\Big(
\widehat{\tau}_ \mathbf{n} =\tau\Big)=1
\end{equation}
and
\begin{equation}\label{conv2}
\lim_{n\rightarrow +\infty}P\Big( \widehat
{s}_  \mathbf{n}  =s\Big)=1.
\end{equation}
Proof of \eqref{conv1}: Since
$\left\{ \widehat{\phi}_{\tau\left(  1\right)  }^{\left( \mathbf{n}\right)
}>\widehat{\phi}_{\tau\left(  2\right)  }^{\left( \mathbf{n}\right)  }%
>\cdots>\widehat{\phi}_{\tau\left(  p\right)  }^{\left(  \mathbf{n}\right)
}   \right\}  \subset\left\{  \widehat{\tau}_\mathbf{n}
=\tau\right\}$, \eqref{conv1} is obtained if we prove  that $\lim_{n\rightarrow+\infty}P\left(    \widehat{\phi}_{\tau\left(  1\right)  }^{\left( \mathbf{n}\right)
}>\widehat{\phi}_{\tau\left(  2\right)  }^{\left( \mathbf{n}\right)  }%
>\cdots>\widehat{\phi}_{\tau\left(  p\right)  }^{\left(  \mathbf{n}\right)
}\right)  =1$ . Note that
\begin{align*}
&  \left\{ \widehat{\phi}_{\tau\left(  1\right)  }^{\left( \mathbf{n}\right)
}>\widehat{\phi}_{\tau\left(  2\right)  }^{\left( \mathbf{n}\right)  }%
>\cdots>\widehat{\phi}_{\tau\left(  p\right)  }^{\left(  \mathbf{n}\right)
}   \right\} 
\\
&  =\left(  \bigcap_{\ell\in E}\bigcap_{i\in F_{\ell}}\left\{  \widehat{\phi
}_{\tau\left(  i\right)  }^{\left(  \mathbf{n}\right)  }>\widehat{\phi}%
_{\tau\left(  i+1\right)  }^{\left( \mathbf{n}\right)  }\right\}  \right)
\cap\left(  \bigcap_{1\leq \ell\leq r-1}\left\{  \widehat{\phi}_{\tau\left(
m_{1}+\cdots+m_{\ell}\right)  }^{\left(  \mathbf{n}\right)  }>\widehat{\phi}%
_{\tau\left(  m_{1}+\cdots+m_{\ell}+1\right)  }^{\left(  \mathbf{n}\right)  }\right\}
\right).
\end{align*}
A similar reasoning than in \cite{nkiet2012} (see the proof of Eq. (5.12), p. 160) leads to the equality
$P\left(  \lim\inf_{n\in\mathbb{N}^{\ast}}\bigcap_{1\leq \ell\leq r-1}\left\{
\widehat{\phi}_{\tau\left(  m_{1}+\cdots+m_{\ell}\right)  }^{\left( \mathbf{n}\right)
}>\widehat{\phi}_{\tau\left(  m_{1}+\cdots+m_{\ell}+1\right)  }^{\left(
\mathbf{n}\right)  }\right\}  \right)  =1$, what implies
\begin{equation}
\lim_{n\rightarrow+\infty}P\left(  \bigcap_{1\leq \ell\leq r-1}\left\{
\widehat{\phi}_{\tau\left(  m_{1}+\cdots+m_{\ell}\right)  }^{\left(  \mathbf{n}\right)
}>\widehat{\phi}_{\tau\left(  m_{1}+\cdots+m_{\ell}+1\right)  }^{\left(
\mathbf{n}\right)  }\right\}  \right)  =1.
\end{equation}
On the other hand, if $E$ is empty, then $ \bigcap_{\ell\in E}\bigcap_{i\in F_{\ell}}\left\{  \widehat{\phi
}_{\tau\left(  i\right)  }^{\left(  \mathbf{n}\right)  }>\widehat{\phi}%
_{\tau\left(  i+1\right)  }^{\left( \mathbf{n}\right)  }\right\}=\Omega$ so $P\Big( \bigcap_{\ell\in E}\bigcap_{i\in F_{\ell}}\left\{  \widehat{\phi
}_{\tau\left(  i\right)  }^{\left(  \mathbf{n}\right)  }>\widehat{\phi}%
_{\tau\left(  i+1\right)  }^{\left( \mathbf{n}\right)  }\right\}\Big)=1$. If  $E\neq\emptyset  $, then for any  $\ell\in E$ and any
$i\in F_{\ell}$ , one has
\begin{align}
n^{d\gamma}\left(  \widehat{\psi}_{\tau\left(  i\right)  }^{\left(  \mathbf{n}\right)
}-\widehat{\psi}_{\tau\left(  i+1\right)  }^{\left( \mathbf{n}\right)  }\right)   &
=n^{d\gamma}\left(  \widehat{\xi}_{K_{\tau\left(  i\right)  }}^{\left(
\mathbf{n}\right)  }-\widehat{\xi}_{K_{\tau\left(  i+1\right)  }}^{\left(  \mathbf{n}\right)
}\right) +f\left(  \tau\left(  i\right)  \right)
-f\left(  \tau\left(  i+1\right)  \right)  .
\end{align}
Then, from Lemma \ref{diffcrit} it follows that   $n^{d\gamma}\left(  \widehat{\psi}_{\tau\left(  i\right)
}^{\left( \mathbf{n}\right)  }-\widehat{\psi}_{\tau\left(  i+1\right)  }^{\left(
\mathbf{n}\right)  }\right)  $ converges in probability to   $f\left(  \tau\left(
i\right)  \right)  -f\left(  \tau\left(  i+1\right)  \right)  $ as $n\rightarrow+\infty$. Then, arguing as in \cite{nkiet2012} (see the end of the proof of Eq. (5.13),  p. 161), we deduce that $
\lim_{n\rightarrow+\infty}P\left(  \widehat{\psi}_{\tau\left(  i\right)
}^{\left( \mathbf{n}\right)  }>\widehat{\psi}_{\tau\left(  i+1\right)  }^{\left(
\mathbf{n}\right)  }\right)  =1$. Since this last property holds  for any  $\ell\in E$ and any
$i\in F_{\ell}$ , we can conclude that $
\lim_{n\rightarrow+\infty}P\left(  \bigcap_{\ell\in E}\bigcap_{i\in F_{\ell}%
}\left\{  \widehat{\psi}_{\tau\left(  i\right)  }^{\left(  n\right)
}>\widehat{\psi}_{\tau\left(  i+1\right)  }^{\left(  n\right)  }\right\}
\right)  =1$ and, finally, that $\lim_{n\rightarrow+\infty}P\left(  \widehat{\phi
}_{\tau\left(  1\right)  }^{\left(  \mathbf{n}\right)  }>\widehat{\phi}%
_{\tau\left(  2\right)  }^{\left(  \mathbf{n}\right)  }>\cdots>\widehat{\phi}%
_{\tau\left(  p\right)  }^{\left(  \mathbf{n}\right)  }\right)  =1$. 

\bigskip

\noindent Proof of \eqref{conv2}:  Clearly,
\[
\left\{  \widehat{\tau}_\mathbf{n}
=\tau\right\}  =\left(  \bigcap_{1\leq
i<s}\left\{  \widehat{\psi}_{s}^{\left(   \textbf{n}\right)  }<\widehat{\psi}%
_{i}^{\left(   \textbf{n}\right)  }\right\}  \right)  \cap\left(  \bigcap_{s<i\leq
p}\left\{  \widehat{\psi}_{s}^{\left(   \textbf{n}\right)  }\leq\widehat{\psi}%
_{i}^{\left(   \textbf{n}\right)  }\right\}  \right)
\]
with the convention $\bigcap_{i\in\emptyset}A_{i}=\Omega$ for any family
$\left(  A_{i}\right)  $ of subsets of $\Omega$. It remains to prove
\begin{equation}\label{conv3}
\lim_{n\rightarrow+\infty}P\left(  \bigcap_{1\leq i<s}\left\{  \widehat{\psi
}_{s}^{\left(   \textbf{n}\right)  }<\widehat{\psi}_{i}^{\left(   \textbf{n}\right)  }\right\}
\right)  =1
\end{equation}
and
\begin{equation}\label{conv4}
\lim_{n\rightarrow+\infty}P\left(  \bigcap_{s<i\leq p}\left\{  \widehat{\psi
}_{s}^{\left(   \textbf{n}\right)  }\leq\widehat{\psi}_{i}^{\left(   \textbf{n}\right)  }\right\}
\right)  =1.
\end{equation}
If $s=1$, then \eqref{conv3}  is obviously obtained from the
equality $\bigcap_{i\in\emptyset}\left\{  \widehat{\psi}_{s}^{\left(
 \textbf{n}\right)  }<\widehat{\psi}_{i}^{\left(   \textbf{n}\right)  }\right\}  =\Omega$. Now
suppose that $s>1$; for any $i\in I$, considering $J_{i}
=\left\{  \tau\left(  \ell\right)  ;\;1\leq \ell\leq i\right\}  $, we have $\left\{  \widehat{\tau}_ \textbf{n} =\tau\right\}  \subset\left\{
\widehat{J}_{i}^{\left(  \textbf{n}\right)  }=J_{i} \right\}
\subset\left\{  \widehat{\xi}_{\widehat{J}_{i}^{\left(   \textbf{n}\right)  }}=\widehat{\xi}_{J_{i}}\right\} $.
Then, from \eqref{conv1}, we deduce that $\widehat{\xi}_{\widehat{J}_{i}^{\left(   \textbf{n}\right)  }}-\widehat{\xi}_{J_{i}}$ converges in probability to $0$ as $n\rightarrow+\infty$, what
implies the convergence in probability of
$\widehat{\xi}_{\widehat{J}_{i}^{\left(   \textbf{n}\right)  }} $ to $ \xi_{J_{i}}$ as  $n\rightarrow+\infty$ since $\widehat{\xi}_{J_{i}}$ almost surely converges  to $ \xi_{J_{i}}$. Further, since
\begin{equation}\label{incl}
  \left\{  \widehat{\tau}^{\left(   \textbf{n}\right)  }=\tau\right\}  \subset\left\{  n^{-d\beta}\Big(g\left(  \widehat{\tau}^{\left(   \textbf{n}\right)  }\left(
i\right)  \right)  -g\left(  \tau\left(  i\right)  \right) \Big) =0\right\}
=\left\{  g\left(  \widehat{\tau}^{\left(   \textbf{n}\right)
}\left(  i\right)  \right)  -g\left(  \tau\left(  i\right)  \right)
 =0\right\}   ,
\end{equation}
it follows that   $n^{-d\beta}\Big(g\left(  \widehat{\tau}^{\left(   \textbf{n}\right)  }\left(
i\right)  \right)  -g\left(  \tau\left(  i\right)  \right) \Big)$    converges
in probability to $ 0$ as $n\rightarrow+\infty$ and, consequently, that  $n^{-d\beta}g\left(  \widehat{\tau}^{\left(   \textbf{n}\right)  }\left(
i\right)  \right)  $ converges in probability to $\ 0$. Hence   $\widehat{\psi}_{i}^{\left( \textbf{n}  \right)  }$
converges in probability to $\ \xi_{J_{i}}$, as 
$n\rightarrow+\infty$. Then, arguing as in \cite{nkiet2012} (see the proof of Eq. (5.15), p. 161), we obtain 
$\lim_{n\rightarrow+\infty}P\left(  \widehat{\psi}_{s}^{\left(  \textbf{n}\right)
}<\widehat{\psi}_{i}^{\left(  \textbf{n}\right)  }\right)  =1$, what gives \eqref{conv3} since this  equality holds for any
\ $i\in\left\{  1,\cdots,s-1\right\}  $. If $s=p$, then \eqref{conv4}  is obviously obtained from the
equality $\bigcap_{i\in\emptyset}\left\{  \widehat{\psi}_{s}^{\left(
\textbf{n}\right)  }\leq\widehat{\psi}_{i}^{\left(  \textbf{n}\right)  }\right\}  =\Omega$. Now
suppose that $s<p$; then for any integer \ $i\in\left\{  s+1,\cdots
,p\right\}  $, we have:
\begin{equation}\label{egalite}
n^{d\beta}\left(  \widehat{\psi}_{i}^{\left(  \textbf{n}\right)  }-\widehat{\psi}%
_{s}^{\left(  \textbf{n}\right)  }\right)  =n^{d\beta}\left(  \widehat{\xi}_{\widehat
{J}_{i}^{\left(  \textbf{n}\right)  }}^{\left(  \textbf{n}\right)  }-\widehat{\xi}_{\widehat
{J}_{s}^{\left(  \textbf{n}\right)  }}^{\left(  \textbf{n}\right)  }\right)  + 
g\left(  \widehat{\tau}^{\left(  \textbf{n}\right)  }\left(  i\right)  \right)
-g\left(  \widehat{\tau}^{\left(  \textbf{n}\right)  }\left(  s\right)  
\right)  .
\end{equation}
First,  for any $k\in\left\{  s,\cdots,p\right\}  $ the  inclusion
$\left\{  \widehat{\tau}_ \textbf{n}   =\tau\right\}  \subset\left\{
\widehat{J}_{k}^{\left(  \textbf{n}\right)  }=J_{k} \right\}
\subset\left\{  n^{d\beta}\left(  \widehat{\xi}_{\widehat{J}_{k}^{\left(
\textbf{n}\right)  }}^{\left(  \textbf{n}\right)  }-\widehat{\xi}_{J_{k} }^{\left(  \textbf{n}\right)  }\right)  =0\right\}$
together with \eqref{conv1}  implies that $n^{d\beta
}\left(  \widehat{\xi}_{\widehat{J}_{k}^{\left(  \textbf{n}\right)  }}^{\left(
 \textbf{n}\right)  }-\widehat{\xi}_{J_{k}}^{\left(  \textbf{n}\right)
}\right)  $ converges in probability to $0$ as $n\rightarrow+\infty$. Further,
since   $I_{1}\subset J_{k} $ (because $k\geq
s$) \ it follows that  $\xi_{J_{k} }=0$, what is equivalent to    $\mathcal{D}_{J_{k}}=0$. Then, from \eqref{dvlpmt}, $n^{d/2}\widehat{\xi}_{J_{k}}^{\left(  \textbf{n}\right)  }$ converges in distribution to
$\left\lVert  \Phi _{J_k}(H)    \right\lVert_\mathcal{H}$ as $n\rightarrow
+\infty$.
Therefore, the equality
$n^{d\beta}\widehat{\xi}_{\widehat{J}_{k}^{(  n)  }}^{(
n)  }=n^{d\beta}(  \widehat{\xi}_{\widehat{J}_{k}^{(  \textbf{n})
}}^{(  \textbf{n})  }-\widehat{\xi}_{J_{k}
}^{(  \textbf{n})  })  +n^{d(\beta-1/2)}( n^{d/2}\widehat{\xi}
_{J_{k}}^{(  \textbf{n})  })$
shows that $n^{d\beta}\widehat{\xi}_{\widehat{J}_{k}^{\left(  \textbf{n}\right)  }%
}^{\left(  \textbf{n}\right)  }$ converges in probability to $0$ as $n\rightarrow
+\infty$. We deduce that
$n^{d\beta}\left(  \widehat{\xi}_{\widehat{J}_{i}^{\left(  \textbf{n}\right)  }%
}^{\left(  \textbf{n}\right)  }-\widehat{\xi}_{\widehat{J}_{s}^{\left(  \textbf{n}\right)  }%
}^{\left(  \textbf{n}\right)  }\right)  $ converges in probability to $0$ as
$n\rightarrow+\infty$. Furthermore,  \eqref{incl} and \eqref{conv1}  show that    $g\left(  \widehat{\tau}^{\left(  \mathbf{n}\right)  }\left(  i\right)  \right)  $
converges in probability to $\ g\left(  \sigma\left(  i\right)  \right)  $ as
$n\rightarrow+\infty$. Consequently, from \eqref{egalite}  it follows that
$n^{d\beta}\left(  \widehat{\psi}_{i}^{\left(  \mathbf{n}\right)  }-\widehat{\psi}%
_{s}^{\left(  \mathbf{n}\right)  }\right)  $ converges in probability to $g\left(
\sigma\left(  i\right)  \right)  -g\left(  \sigma\left(  s\right)  \right)  $
as $n\rightarrow+\infty$. Then, arguing as in \cite{nkiet2012} (see the proof of Eq. (5.16), p. 162), we obtain  $\lim_{n\rightarrow+\infty}P\left(   \widehat{\psi}_{s}^{\left(
\mathbf{n}\right)  }\leq\widehat{\psi}_{i}^{\left( \mathbf{n}\right)  }  \right)
=1$, what yields \eqref{conv4} since  this latter equality holds for any $i\in\left\{
s+1,\cdots,p\right\}  .$





\end{document}